\input amstex
\magnification=\magstep{1.5}
\baselineskip 22pt
\documentstyle{amsppt}
\topmatter
\title 
Calabi-Yau threefolds of quotient type
\endtitle
\author 
Keiji Oguiso and Jun Sakurai
\endauthor
\affil 
Mathematisches Institut Universit\"at Essen, Department of Mathematical 
Sciences University of Tokyo
\endaffil
\address
D-45117 Essen Germany, 153-8914 Komaba Meguro Tokyo Japan
\endaddress
\email
mat9g0\@spi.power.uni-essen.de
\endemail
\dedicatory
Dedicated to Professor V. A. Iskovskih on the occasion of his sixtieth birthday
\enddedicatory
\abstract 
By a Calabi-Yau threefold we mean a minimal complex projective 
threefold $X$ such that $\Cal O_{X}(K_{X}) \simeq \Cal O_{X}$ and 
$h^{1}(\Cal O_{X}) = 0$. This paper consists of four parts. 
In Section 1 we formulate an equivariant version of Torelli Theorem 
of K3 surface with a finite group action and 
deduce some more geometrical consequences.  
In Section 2, we classify Calabi-Yau threefolds with infinite 
fundamental group by means of their minimal splitting coverings introduced 
by Beauville, and deduce as its Corollary that the nef 
cone is a rational simplicial cone and any rational nef divisor 
is semi-ample provided that $c_{2}(X) \equiv 0$ on 
$\text{Pic}(X)_{\Bbb R}$. We also derive a 
sufficient condition for $\pi_{1}(X)$ to be finite in terms of 
the Picard number in an optimal form. 
In Section 3, we give a fairly concrete structure Theorem concerning 
$c_{2}$-contractions of Calabi-Yau threefolds as a generalisation 
and also a correction of our 
earlier works for simply connected ones. 
In Section 4, applying the results in these three 
sections together with Kawamata's finiteness result of 
the relatively minimal models of a Calabi-Yau fiber space, 
we show the finiteness of the isomorphism classes of 
$c_{2}$-contractions of each Calabi-Yau threefold. 
As a special case, we find the finiteness of 
abelian pencil structures on each $X$ up to $\text{Aut}(X)$.   
\endabstract
\toc 
\subhead
\S 0. Introduction
\endsubhead
\subhead
\S 1. An equivariant Torelli Theorem for K3 surfaces with 
finite group action and its applications
\endsubhead
\subhead
\S 2. Calabi-Yau threefolds with infinite fundamental group 
\endsubhead
\subhead
\S 3. Classification of $c_{2}$-contractions of Calabi-Yau threefolds
\endsubhead
\subhead
\S 4. Finiteness of $c_{2}$-contractions of a Calabi-Yau threefold
\endsubhead
\endtoc 
\endtopmatter
\document
\head
{\S 0. Introduction}
\endhead 
In the light of the minimal model theory, we define {\it a Calabi-Yau threefold} to be a $\Bbb Q$-factorial terminal projective threefold $X$ defined over 
$\Bbb C$ such that $\Cal O_{X}(K_{X}) \simeq 
\Cal O_{X}$ and $h^{1}(\Cal O_{X}) = 0$, and regard the second Chern class 
$c_{2}(X)$ as a linear form on $\text{Pic}(X)_{\Bbb R}$ 
through the intersection 
pairing, where $c_{2}(X)$ for a singular $X$ is defined as $c_{2}(X) := 
\nu_{*}(c_{2}(\tilde{X}))$ via a resolution $\nu : \tilde{X} \rightarrow X$ 
and is known to be well-defined (see for example [Og1, Lemma (1.4)]). 
\par
\vskip 4pt 
However, as is pointed out by several authors, this preferable definition 
of Calabi-Yau threefold has an inevitable defect: Those Calabi-Yau 
threefolds, such as Igusa's example ([Ig, Page 678], [Ue, Example 16.16]), 
that are given as an \'etale quotient of an abelian threefold are then 
included in our category. We call them {\it of Type A}. Indeed, their 
pathological 
nature sometimes prevents us from studying Calabi-Yau threefolds uniformly. 
For instance,
\roster
\item 
there are no rational curves on Calabi-Yau threefolds of Type A, while it is 
expected, and has been already checked in some extent, that most of Calabi-Yau 
threefolds contain rational curves (see [Wi1], [HW] and [EJS]);
\item 
$c_{2}(X) = 0$ for such $X$ but $c_{2}(X) \not= 0$ for others.
\endroster
Here, for the last statement, we recall the following result due to 
S. Kobayashi in smooth case and Shephard-Barron and Wilson in the general 
case:
\proclaim{Theorem ([Kb, Chap. IV, Corollary (4.15)], [SBW, Corollary])} 
Let $X$ be a Calabi-Yau threefold. Then, $X$ is of 
Type A if and only if $c_{2}(X) = 0$. \qed
\endproclaim 
One of the main purposes of this paper is to compensate for this defect by 
revealing explicit geometric structures of Calabi-Yau threefolds of Type A. 
It will turn out that they are remarkably few so 
that one can in principle handle them separately in case. Our result is:

\proclaim{Theorem (0.1)}  
Let $X$ be a Calabi-Yau threefold of Type A.  
Then, 
\flushpar 
(I) $X = A/G$, where $A$ is an abelian threefold and $G$ is its finite 
automorphism group acting freely on $A$ such that either one of 
the following (1) or (2) is satisfied: 
\roster
\item 
$G = \langle a \rangle \oplus \langle b \rangle \simeq C_{2}^{\oplus 2}$ and,
\par
$a_{0}  = 
\left(\matrix 1 & 0 & 0 \\ 0 & -1 & 0 \\ 0 & 0 & -1 \endmatrix \right)$ 
and 
$b_{0} = \left(\matrix -1 & 0 & 0 \\ 0 & 1 & 0 \\ 0 & 0 & -1\endmatrix \right)$
\par    
\item 
$G = \langle a, b \vert a^{4} = b^{2} = 1, bab = a^{-1} \rangle \simeq D_{8}$ 
and,
\par
$a_{0} = 
\left(\matrix 1 & 0 & 0 \\ 0 & 0 & -1 \\ 0 & 1 & 0 \endmatrix \right)$ 
and 
$b_{0} 
= \left(\matrix -1 & 0 & 0 \\ 0 & 1 & 0 \\ 0 & 0 & -1\endmatrix \right)$, 
\endroster
where $a_{0}$ and $b_{0}$ stand for the Lie parts of $a$ and $b$ 
respectively and the matrix representation is the one given under appropriate 
global coordinates of $A$.  
\flushpar
(II) In the first case $\rho(X) = h^{1}(T_{X}) = 3$ and in the second case 
$\rho(X) = h^{1}(T_{X}) = 2$. In particular, a Calabi-Yau threefold with 
Picard number $\rho \not= 2, 3$ is not of Type A. 
\flushpar
(III) Both cases actually occur (see (2.17) and (2.18) for 
explicit examples). 
\flushpar
(IV) In each case, the nef cone $\overline{\Cal A}(X)$ is a rational, 
simplicial cone and every rational nef divisor on $X$ is semi-ample. 
In particular, $X$ admits only finitely many contractions. 
\endproclaim 

We call a pair $(A, G)$ which falls into the cases (1) and (2) 
{\it Igusa's pair} and {\it refined Igusa's pair} respectively. 
\par
\vskip 4pt 
This Theorem is shown in the subsection (2A) based on the notion of 
the minimal splitting covering introduced by Beauville ([Be2, Section 3], 
see also 
(2.1)), which nicely reduces a greater part of the proof to the problem 
of representations of a special kind of groups $G$ called {\it pre-Calabi-Yau 
groups of Type A }((2.2)). Then, one of the important steps is to restrict the 
possible orders of such $G$. For this purpose, we apply the 
Burnside-Hall Theorem concerning commutative subgroups of $p$-groups 
([Su, Page 90], see also (2.7)). 
\par
\vskip 4pt 
Let us add a few remarks about the result. 
First, even compared with the range $1 \leq \rho \leq 9$ of the Picard 
numbers of abelian threefolds, the range $\{2, 3 \}$ of $\rho(X)$ is quite 
narrow. 
Secondly, our statements (II) and (III) show that there certainly exist 
smooth Calabi-Yau threefolds containing no rational curves if $\rho = 2$ and 
$3$, but, on the other hand, suggest some hope to ask the following: 

\example{Question (cf. [Wi1], [MS])} Does every Calabi-Yau 
threefold of Picard number $\rho \not= 2, 3$ contain rational curves? 
\qed
\endexample 
The third remark concerns the statement (IV). Recall that any smooth 
anti-canonical member $X \in \vert -K_{V} \vert$ of a smooth Fano fourfold 
$V$ is a simply connected Calabi-Yau threefold.  
In addition, such an $X$ always satisfies $c_{2}(X) > 0$ on 
$\overline{\Cal A}(X) - \{0\}$ ([OP, Main Theorem 2]). 
So, the statement (IV) can be regarded as an extreme counterpart of the 
following Theorem due to Koll\'ar: 

\proclaim{Theorem ([Bo, Appendix])} Let $X$ be a smooth member of 
$\vert -K_{V} \vert$ of a smooth Fano fourfold $V$ and $\iota : X \rightarrow 
V$ the natural inclusion. Then $\iota^{*} : \overline{\Cal A}(V) \rightarrow 
\overline{\Cal A}(X)$ is an isomorphism. In particular, 
the nef cone $\overline{\Cal A}(X)$ is a rational polyhedral cone, 
every rational nef divisor on $X$ is semi-ample, and therefore, $X$ 
admits only finitely many contractions. \qed
\endproclaim 

Refer also to [Wi1, Page 146, Claim], [Wi2, 4] and [Og1, Theorem (2.1), 
Proposition (2.7)] for related results about semi-ampleness of nef 
divisors.
\par
The last remark is concerned with fundamental group. According to the 
Bogomolov decomposition Theorem [Be1] and its generalisation due to 
Yoshinori Namikawa and Steenbrink [NS, Corollary (1.4)], there is 
one more class of Calabi-Yau threefolds with infinite fundamental group, 
namely, the class consisting of those Calabi-Yau threefolds that are given as 
an \'etale quotient of $\text{(K3 surface)}\, \times\, 
\text{(elliptic curve)}$. We call 
them Calabi-Yau threefolds {\it of Type K} and study them in subsection 
(2K) in some extent. (See Theorem (2.23) for the statement.) As an 
application of (0.1) and (2.23), we obtain the following criterion for 
$\pi_{1}(X)$ 
to be finite in terms of the Picard number: 

\proclaim{Corollary (0.2)} 
Let $X$ be a Calabi-Yau threefold. 
Then, $\pi_{1}(X)$ is finite if $\rho(X) = 1, 6, 8, 9, 10$ or 
$\rho(X) \geq 12$. 
Moreover, this is also optimal in the sense that for each 
$\rho \in \Bbb N - \{1, 6, 8, 9, 10, n \geq 12\} (= \{2, 3, 4, 5, 7, 11\})$, 
there exists a Calabi-Yau threefold $X$ such that $\pi_{1}(X)$ is infinite 
and $\rho(X) = \rho$. In particular, the fundamental group of 
a Calabi-Yau threefold whose Picard number one is always finite.  
\endproclaim 

The last statement in (0.2) is also obtained by Amerik, 
Rovinsky and Van de Ven [ARV, Proposition (3.1)], which the first author 
was kindly informed by Amerik after his talk on this subject. 
\par
\vskip 4pt 
So far, we have concerned special kinds of Calabi-Yau threefolds called 
of Type A and of Type K.  Another particular interest of this paper, which 
turns out to be related to our first problem, is the role of the second 
Chern class $c_{2}(X)$ in the geometry of contractions of $X$. 
Here, the term {\it contraction} means a surjective morphism onto a normal, 
projective variety with connected fibers, and therefore, consists of the two 
cases, that is, the fiber space case and the birational contraction case. 
Let $\varphi : X \rightarrow W$ be a contraction. Then, there is a nef 
divisor $D$ on $X$ such that $\varphi = \Phi_{D}$, where $\Phi_{D}$ stands 
for the morphism associated with the complete linear system 
$\vert D \vert$. Therefore, we may relate $\varphi$ with $c_{2}(X)$ via the 
intersection number $(c_{2}(X).D)$. Although 
the value $(c_{2}(X).D)$ itself is not well-defined for $\varphi$, 
it does not depend on the choice of $D$ 
such that $\varphi = \Phi_{D}$ whether $(c_{2}.D) = 0$ or not. 
This is due to the pseudo-effectivity of $c_{2}(X)$ ([Mi]). We call 
$\varphi$ a $c_{2}$-{\it contraction} if $(c_{2}(X).D) = 0$. 
For example, a pencil $\varphi : X \rightarrow \Bbb P^{1}$ is 
$c_{2}$-contraction if and only if the general fiber of $\varphi$ is 
an abelian surface.  
\par
\vskip 4pt
Our first task in this direction is to enlarge our earlier classification of 
$c_{2}$-contractions in simply connected case ([Og 2, 3, 4]) to the one 
in the general case as in (0.3) below. For the statement, we recall the 
following 
pairs of 
an abelian threefold and its specific Gorenstein automorphism: 
the pair $(A_{3}, g_{3})$, where 
$A_{3}$ is the 
product threefold $A_{3} := E_{\zeta_{3}}^{3}$ of the elliptic curve 
$E_{\zeta_{3}}$ of period $\zeta_{3} = \text{exp}(2\pi i/3)$ and $g_{3}$ 
is its automorphism $\text{diag}(\zeta_{3}, \zeta_{3}, \zeta_{3})$; and the 
pair $(A_{7}, g_{7})$, where $A_{7}$ is the Jacobian threefold of the Klein 
quartic curve $C = (x_{0}x_{1}^{3} + x_{1}x_{2}^{3} + x_{2}x_{0}^{3} = 0) 
\subset 
\Bbb P^{2}$ and $g_{7}$ is the automorphism of $A_{7}$ induced by the 
automorphism of $C$ given by $[x_{0} : x_{1} : x_{2}] \mapsto 
[\zeta_{7}x_{0} : \zeta_{7}^{2}x_{1} : \zeta_{7}^{4}x_{2}]$. 
We call 
$(A_{3}, g_{3})$ {\it the Calabi pair} and $(A_{7}, g_{7})$ 
{\it the Klein pair}.  

\proclaim{Theorem (0.3)} Let $X$ be a Calabi-Yau threefold. Assume that 
$X$ admits a $c_{2}$-contraction $\varphi : X \rightarrow W$ such that 
$\text{dim}(W) \geq 2$. Then, $X$ is smooth and is birational to either 
one of the following: 
\roster 
\item 
a crepant resolution of a Gorenstein quotient $(S \times E)/G$ of the product 
of a normal K3 surface $S$ and an elliptic curve $E$, where by a normal K3 
surface we mean a normal surface whose minimal resolution is a K3 surface; 
\item 
the crepant resolution of $A/G$, where $(A,G)$ is either the Calabi pair, its 
modification, the Klein pair, Igusa's pair or refined Igusa's pair.
\endroster 
(See (3.3), (3.4), (3.6) and (3.7) for more precise statement and 
structures.) 
\endproclaim
The main idea of the proof of Theorem (0.3) is to modify $\varphi : 
X \rightarrow W$ 
toward one of the threefolds described in (0.3) by 
taking appropriate coverings, Stein factorisations, or by running log 
mimimal model program, as intrinsically as possible in order 
to inherite group actions biregularly. 
This idea itself is same as the one 
in simply connected case and, indeed, most part of the proof toward 
the case (2) can be done by 
a combination of (0.1) and more or less obvious minor modification 
of [Og2, 3].    
However, the first author should confess that his argument of 
[Og4, Section 3] concerning lifting of certain group actions {\it contains a 
gap}, which he noticed around August 1999, 
and the argument [Og4, Sections 3, 4] toward the case (1) is not available. 
(Claim (3.4) in [Og4] is false and 
the right statement is that $\nu$ in (3.4) is at best the normalisation 
{\it in a certain finite field extension}. Therefore,  
$\nu$ is not intrinsic so that the lifting argument there and the 
argument after that seem to be broken.)  Unlike the argument 
there, which is based on the theory of quasi-product 
threefold, our new 
idea here is to apply again the notion of minimal splitting covering of 
Beauville, especially, its uniqueness property, after the same reduction 
as in [Og4, Section 2]. Fortunately, this argument also recovers the main 
result [Og4] as its own form and even simplifies the proof. This will be 
done in Section 3, especially in (3.7).
\par
\vskip 4pt
The final aim of this paper is to show the following: 
\proclaim{Theorem (0.4)} 
Each Calabi-Yau threefold $X$ admits only finitely many different 
$c_{2}$-contractions up to isomorphism. In particular, $X$ admits 
only finitely many different abelian pencil structures up to $\text{Aut}(X)$. 
\endproclaim 
\par
\vskip 4pt 
This result is particularly motivated by the following work on the finiteness 
of fiber spaces coming from "the opposite side" of the cone:  
\proclaim{Theorem [OP, Main Theorem 1 and Remark in Section 3]}
Let $X$ be a Calabi-Yau threefold and let $H$ be an ample divisor on $X$ and 
$\epsilon > 0$ a positive real number. 
Set $\overline{\Cal A}_{\epsilon}(X) := \{x \in \overline{\Cal A}(X) \vert 
(c_{2}(X).x) \geq \epsilon(H^{2}.x) \}$. 
\roster
\item 
Assume that $c_{2}(X) > 0$ on $\overline{\Cal A}(X) - \{0\}$. 
Then, $X$ admits only finitely many different fibrations. 
\item 
More generally, the cardinality of the fibrations $\varphi : X \rightarrow W$ 
such that 
$\varphi^{*}\overline{\Cal A}(W) \subset \overline{\Cal A}_{\epsilon}(X)$ 
is finite.  \qed
\endroster
\endproclaim 
Both Theorem (0.4) and Theorem [OP] are related positively to the Cone 
Conjecture posed by D. Morrison [MD]. However, these two Theorems are 
completely different both in nature and in proof. Proof of Theorem [OP] 
is based on the boundedeness results of log surfaces due to Alexeev [Al] and 
the compactness of the domain $\{x \in \overline{\Cal A}_{\epsilon}(X) 
\vert (c_{2}(X).x) \leq B\}$. Compactness, in particular, implies the 
finiteness of the lattice points in the domain. Main idea in [OP] is 
to reduce the problem to this finiteness by applying the boundedness.  
Therefore, the result claims the finiteness of fiber spaces in question 
themselves. (See [OP] for details.) However, this compactness reduction 
does not work any more for $c_{2}$-contractions. 
In addition, there actually exists a Calabi-Yau threefold which admits 
infinitely many different abelian pencils ([Og1, Section 4]). 
Therefore, contrary to [OP], the finiteness of $c_{2}$-contractions 
themselves are false in general and it should be the core of (0.4) 
that we modulo them up by automorphisms. 
Outline of proof of (0.4) is as follows: 
We first take {\it the maximal} $c_{2}$-{\it contraction of} $X$ ((4.1)) 
and denote this by $\varphi_{0} : X \rightarrow W_{0}$. This $\varphi_{0}$ 
has a property that any $c_{2}$-contraction factors through $\varphi_{0}$. 
Then we divide into cases according to the structure of $\varphi_{0}$. 
The essential case is the case where $\varphi_{0}$ falls into the case (1) 
of (0.3). In this case, we devide our problem into two parts: 
finiteness of contractions of $W_{0}$ up to isomorphisms; and, 
lifting of automorphisms of $W_{0}$ to $\varphi_{0} : X \rightarrow W_{0}$. 
We apply Kawamata's finiteness result of relatively minimal models of 
a Calabi-Yau fiber space [Kaw5, Theorem (3.6)] to the second part 
and an equivariant Torelli Theorem for pair $(S, G)$ 
of K3 surface and its finite automorphism group ((1.8), (1.10)) 
to the first part. The reason why the finiteness of minimal models is 
needed is because minimal models of $X$ are no more unique. 
The role of our equivariant Torelli Theorem may be more or less apparent 
by the fact that $W_{0}$ is birational to $S/G$ for some 
$(S, G)$. Full verification of (0.4) will be given in Section 4. 
It might be also worth while noticing here that, 
in order to examine automorphisms toward finiteness, 
Kawamata [Kaw5] extended 
Torelli Theorem to the one over non-closed field and applied ``vertically'' 
to his finiteness result,  
while we extended it to the one with finite group action and 
applied ``horizontally'' to our finiteness result. 
\par
\vskip 4pt
We formulate our equivariant Torelli Theorem in Section 1. 
Besides the present application (see also (2.23)(IV)), this Theorem 
has been also applied to study finite automorphism groups of K3 surfaces 
by [OZ1, 2, 3].  

\subhead
Acknowledgement 
\endsubhead
\par 
\vskip 4pt
The present work was initiated during 
the first named author's stay at the Johns Hopkins University in April-June 
1996 and the University of Warwick in July 1996, the main idea of Section 4 
for simply-connected Calabi-Yau threefolds (but based on [Og4]) has been 
found during his stay at the Institut Mittag-Leffler in May 1997,
 an important example (2.33) for Corollary 
(0.2), has been found during the first named author's stay in the 
Universit\"at 
Bayreuth July 1999, and the present form of this work included the final proof 
of the results in Section 3  has been found and completed during his stay in 
Universit\"at Essen 1999 under the financial support by the 
Alexander-Humboldt fundation.  
The first named author would like to express gratitude to 
Professors Y. Kawamata and V. V. Shokurov; 
Professors Y. Miyaoka and M. Reid; 
Professors E. Ellingsrud, D. Laksov, A. Str$\phi$mme and A. Thorup; 
Professor T. Peternell and 
Professors H. Esnault and E. Viehweg 
for their invitations and their warm hospitalities. 
Deep appreciation goes to the Johns Hopkins University, JSPS programmes, 
the University of Warwick, the Institut Mittag-Leffler, Universit\"at 
Bayreuth and the Alexander-Humboldt fundation for their financial supports, 
and Universit\"at Essen for providing a pleasant environment to research. 
A greater part of Section 2 has grown out of the second named author's master 
thesis at University of Tokyo 1998 under the first named author's 
instruction. Last but not least at all, both authors would like to express 
their thanks to Professor Y. Kawamata for his warm encouragement during 
their preparation and to Professors E. Amerik, F. Campana, V. V. Nikulin, 
D. Q. Zhang and Q. Zhang for their comments and encouragement. 
\par 
\vskip 4pt
Thoughout this paper, in addition to the notation and terminology in 
[Ha], [KMM] and in Introduction, we employ the following:
\subhead
{\bf Notation, Terminology and Convention} 
\endsubhead  
\definition{(0.1)} Every variety in this paper is assumed to be normal and 
projective defined over $\Bbb C$ unless stated otherwise. The open convex cone 
generated by the ample classes in $N^{1}(X) := 
(\{\text{Cartier\, divisors}\}/\equiv) 
\otimes \Bbb R$ is called the ample 
cone and is denoted by $\Cal A(X)$. Its closure $\overline{\Cal A}(X)$ is 
called {\it the nef cone}. A $\Bbb Q$-Cartier divisor $D$ on $X$ is said to be 
{\it semi-ample} 
if there exists a positive integer $m$ such that $\vert mD \vert$ is free. 
\enddefinition
\definition{(0.2)} Two contractions $\varphi : X \rightarrow W$ and 
$\varphi' : X' \rightarrow
W'$ are said to be {\it isomorphic} if there exist isomorphisms 
$F : X \rightarrow X'$ and $f : W \rightarrow W'$ such that 
$\varphi' \circ F = f \circ \varphi$: 
$$\CD 
X @> F >> X' \\
@V \varphi VV @VV \varphi' V \\
W @>> f > W'.
\endCD$$ 
Two contractions of $X$, 
$\varphi : X \rightarrow W$ and $\varphi' : X \rightarrow W'$ are said to 
be {\it identically isomorphic} if there exists an isomorphism 
$f : W \rightarrow W'$ such that $\varphi' = f \circ \varphi$. 
Identically isomorphic contractions should be considered to be the same.  
It is important to distinguish these two notions, 
isomorphic and identically isomorphic, especially in 
the case where $X = X'$. For example, two natural projections 
$p_{i} : X \times X \rightarrow X$ are clearly isomorphic but never 
identically isomorphic. Their difference in terms of the 
nef cone is as follows: Two contractions $\varphi : X \rightarrow W$ and 
$\varphi' : X \rightarrow W'$ are isomorphic if and only if there exists 
an autmorphism $F \in \text{Aut}(X)$ such that $F^{*}\varphi^{*}\overline{\Cal 
A}(W) = (\varphi')^{*}\overline{\Cal A}(W')$, while these two are 
identically isomorphic if and only if $\varphi^{*}\overline{\Cal A}(W) 
= (\varphi')^{*}\overline{\Cal A}(W')$. 
\enddefinition 
\definition{(0.3)} Let $X$ be a Gorenstein variety such that 
$\Cal O_{X}(K_{X}) \simeq \Cal O_{X}$. We denote by $\omega_{X}$ a generater 
of 
$H^{0}(\Cal O_{X}(K_{X}))$. A finite automorphism group $G \subset 
\text{Aut}(X)$ is called {\it Gorenstein} if $g^{*}\omega_{X} = \omega_{X}$ 
for each $g \in G$. 
\enddefinition 
\definition{(0.4)} Throughout this paper, we often need to examine an object 
with a faithful group action $G \curvearrowright S$. We distinguish several 
notions concerning ``fixed loci'' by using the following different symbols:
\flushpar
$S^{g} := \{s \in S \vert g(s) = s\}$ for $g \in G$;
\flushpar 
$S^{G} := \cap_{g \in G} S^{g}$, the set of points which are fixed by 
{\it all} the 
elements of $G$; 
\flushpar
$S^{[G]} := \cup_{g \in G -\{1\}} S^{g}$, the set of points which are fixed by 
{\it some} non-trivial element of $G$. An action $G \curvearrowright S$ is 
said to be fixed point free if $S^{[G]} = \emptyset$.
\enddefinition 
\definition{(0.5)} Let $G$ be a finite group and $X$ a variety with a group 
action $\rho_{X} : G \rightarrow \text{Aut}(X)$. A contraction $\varphi : X 
\rightarrow W$ is said to be 
$G$-{\it stable} if there exists 
a representation $\rho_{W} : G \rightarrow 
\text{Aut}(W)$ which satisfies $\varphi \circ \rho_{X}(g) = \rho_{W}(g) \circ 
\varphi$. Note that the representation $\rho_{W}$ is uniquely determined 
by $\rho_{X}$ and $\varphi$. The $G$-stability of the contraction is also 
equivalent to the existence of a line bundle $D \in \text{Pic}(X)^{G}$ 
such that $\varphi = \Phi_{D}$. 
Two $G$-stable contractions $\varphi : X \rightarrow W$ and 
$\varphi' : X' \rightarrow W'$ are said to be $G$-{\it equivariantly 
isomorphic} 
if there exist isomorphisms 
$F : X \rightarrow X'$ 
and $f : W \rightarrow W'$ such that 
$$\varphi' \circ F = f \circ \varphi,\, F \circ \rho_{X}(g) = 
\rho_{X'}(g) \circ F,\, f \circ \rho_{W}(g) = \rho_{W'}(g) \circ f.$$ 
\enddefinition 
\definition{(0.6)} $\zeta_{n} := \text{exp}(2\pi \sqrt{-1}/n)$, the primitive 
$n$-th root 
of unity in $\Bbb C$.
\enddefinition
\definition{(0.7)} We denote some specific groups appearing in this paper 
by the following fairly standard symbols:
\flushpar
$C_{n} := \langle a \vert a^{n} = 1 \rangle$, the cyclic group of order $n$;
\flushpar
$D_{2n} := \langle a, b \vert a^{n} = b^{2} = 1, bab = a^{-1} \rangle \simeq 
C_{n} \rtimes C_{2}$, the dihedral group of order $2n$; 
\flushpar
$Q_{4n} := \langle a, b \vert a^{2n} = 1, a^{n} = b^{2}, b^{-1}ab = a^{-1} 
\rangle$, 
the binary dihedral group of order $4n$; 
\flushpar
$S_{n} := \text{Aut}_{\text{set}}(\{1,2, \cdots, n\})$, the $n$-th symmetric 
group;
\flushpar
$A_{n} := \text{Ker}(\text{sgn} : S_{n} \rightarrow \{\pm 1\})$, the $n$-th 
alternative group.
\enddefinition
\definition{(0.8)} Let $A := \Bbb C^{d}/\Lambda$ be a $d$-dimensional 
complex torus. 
By abuse of language, we call global coordinates $(z_{1}, z_{2}, ..., z_{d})$ 
of $\Bbb C^{n}$ {\it global coordinates of} $A$ if they are obtained by an 
affine 
transformation of the natural global coordinates of $\Bbb C^{d}$ given by the 
$i$-th projections. When the origin $0$ of $A$ is specified and $A$ is 
regarded as a group variety, we identify $A$ with its translation group 
in the natural 
manner and denote by $(A)_{n}$ the group of $n$-torsion points and by 
$t_{*}$ the translation given by $* \in A$. Under this identification, 
we have $\text{Aut}(A) = A \rtimes \text{Aut}_{\text{Lie}}(A)$, 
where $\text{Aut}_{\text{Lie}}(A)$ is the subgroup consisting 
of the elements $g$ such that $g(0) = 0$. We often call the second 
factor of $h \in \text{Aut}(A)$ under this decomposition, 
{\it the Lie part of} 
$h$ and denote it by $h_{0}$. We also 
denote by $E_{\zeta}$ the elliptic curve whose period is $\zeta$ in 
the upper half plane. 
\enddefinition
\definition{(0.9)} In this paper, we often regard group actions on varieties 
as the so-called co-action through their coordinates. The advantage of this 
convention is that we may then describe its action on cohomology as if 
it were covariant, namely, $(ab)^{*} = a^{*}b^{*}$.   
\enddefinition
\definition{(0.10)} We abbreviate by $S_{K}$ the scalar extension 
$S \otimes_{\Bbb Z} K$ of the space $S$. 
\enddefinition

\head
{\S 1. An equivariant Torelli Theorem for K3 surfaces 
with finite group action and its applications}
\endhead 
\par 
\vskip 4pt
Let $X$ be a K3 surface and $G \subset \text{Aut}(X)$ 
a finite automorphism group. Throughout this section this pair 
$(X, G)$ is fixed. 
The aim of this section is to formulate an equivariant version of 
the Torelli Theorem (1.8) which describes the automorphisms of $X$ which 
commute with $G$ in terms of their actions on cohomology, and 
apply this to get more geometrical consequences (1.9)-(1.11). The core of 
the formulation is to define the $G$-equivariant reflection group (1.6). 
 
\definition{(1.1)} As usual, we consider the second cohomology group 
$H^{2}(X, \Bbb Z)$ as a lattice 
by the non-degenerate symmetric bilinear form $(*.*)$ induced by 
the cup product. We denote by $O(H^{2}(X, \Bbb Z))$ 
the orthogonal group of $H^{2}(X, \Bbb Z)$. 
\enddefinition 
\definition{(1.2)} For the convenience of the formulation, we introduce the 
following notation which in principle follows the rule that 
$U$ denotes the $G$-invariant 
part of the abstract one $U'$. (The reason behind this usage of notation is 
the fact that $G$-invariant part plays more important roles 
in our formulation.) Another rule is 
that the symbol $W^{+}$ indicates a subgroup of a group $W$:  
\flushpar 
$S' := H^{1,1}(X) \cap H^{2}(X, \Bbb Z)$, the N\'eron-Severi lattice of $X$;
\flushpar 
$\Cal N' := \{ [E] \in S' \vert E \subset X, E \simeq \Bbb P^{1}\}$, 
the set of nodal classes; 
\flushpar
$S := (S')^{G} (= \{x \in S' \vert g^{*}(x) = x \, \text{for all}\,  
g \in G\})$;
\flushpar
$T := S^{\perp} = \{x \in H^{2}(X, \Bbb Z) \vert (x.y) = 0 \, 
\text{for all}\, y \in S\}$, the orthogonal lattice of $S$ in 
$H^{2}(X, \Bbb Z)$ (Note that $T$ contains the transcendental lattice of $X$.);
\flushpar
$S^{*} := \text{Hom}(S, \Bbb Z)$, which we always regard as an 
overlattice of $S$, $S \subset S^{*} 
\subset S_{\Bbb Q}$, via the non-degenerate pairing $(*.*)\vert S$; 
\flushpar
$(\Cal C')^{\circ} :=$ the positive cone of $X$, that is, the connected 
component of the space 
$\{x \in (S')_{\Bbb R} \vert (x.x) > 0\}$ containing 
the ample classes;
\flushpar 
$\Cal C' :=$ the union of $(\Cal C')^{\circ}$ and all $\Bbb Q$-rational 
rays in the boundary $\partial (\Cal C')^{\circ}$ of 
$(\Cal C')^{\circ} (\subset (S')_{\Bbb R})$;
\flushpar
$\Cal C := (\Cal C')^{G} = \Cal C' \cap S_{\Bbb R}$; 
\flushpar
$\Cal A' :=$ the intersection of the nef cone $\overline{\Cal A}(X)$ and 
$\Cal C'$; 
\flushpar 
$\Cal A := (\Cal A')^{G} = \Cal A' \cap S_{\Bbb R}$; 
\flushpar
$Q := \{f \in \text{Aut}(X) \vert f \circ g = g \circ f \, \text{for all}\,  
g \in G \}$; 
\flushpar
$O(S) :=$ the orthogonal group of the lattice $S$ preserving $\Cal C$;
\flushpar 
$O(S)^{+} :=$ the subgroup of $O(S)$ consisting of the elements of 
the form $\tau \vert S$, where $\tau$ is a Hodge isometry 
of $H^{2}(X, \Bbb Z)$ such that $\tau g^{*} = g^{*} \tau$ 
for all $g \in G$ (Note that by the last condition, such $\tau$ always 
satisfy $\tau(S) = S$); 
\flushpar
$P(S) := \{\sigma \in O(S) \vert \sigma(\Cal A) = \Cal A\}$;
\flushpar
$P(S)^{+} := \{\sigma \in O(S) \vert \sigma = f^{*}\vert S \, 
\text{for some}\, f \in Q\}$.
\enddefinition
\par
\vskip 4pt 
We should keep in mind the following easy facts and relations: 
\proclaim{Lemma (1.3)}
\roster
\item 
$S$ is an even hyperbolic lattice if $\text{rank}\, S \geq 2$. 
\item
The interior $\Cal A^{\circ}$ of $\Cal A$ 
consists of the $G$-invariant ample classes of $X$ and is non-empty. 
Moreover, $\Cal A^{\circ} = (\Cal A')^{\circ} \cap S_{\Bbb R}$.
\item 
$f^{*}(S) = S$ and $f^{*} \vert S \in P(S)^{+}$ for all $f \in Q$. 
In other words, $P(S)^{+}$ is the image of the homomorphism 
$Q \rightarrow O(S)$ given by $f \mapsto f^{*}\vert S$. 
In particular, $\sigma(\Cal A) = \Cal A$ for each $\sigma \in P(S)^{+}$.
\item 
Set $O(S)^{++} := \{\sigma \in O(S)\, \vert\,  
\sigma \vert (S^{*}/S) = id\}$. Then, 
$O(S)^{++} \subset O(S)^{+} \subset O(S)$ and each of these 
inclusions is of finite index.
\endroster
\endproclaim
\demo{Proof} The assertions (1) and (2) are clear. Note that 
a Hodge isometry $\sigma$ of $H^{2}(X, \Bbb Z)$ satisfies $\sigma(S) 
= S$ if $\sigma g^{*} = g^{*} \sigma$ for all $g \in G$. This 
implies the assertion (3). We show the assertion (4). 
Since $O(S)^{++} = \text{Ker}(O(S) \rightarrow \text{Aut}(S^{*}/S))$ 
and $\text{Aut}(S^{*}/S)$ is a finite group, $O(S)^{++} \subset O(S)$ 
is of finite index. It remains to show that 
$O(S)^{++} \subset O(S)^{+}$. 
Let $\sigma$ be an element of $O(S)^{++}$. 
Since $\sigma \vert (S^{*}/S) = id$, the pair $(\sigma, id) \in 
O(S) \times O(T)$ can be extended to an element 
$\tau \in O(H^{2}(X, \Bbb Z))$. This $\tau$ is a Hodge isometry, because 
$[\omega_{X}] \in T_{\Bbb C}$ and $\tau \vert T = id$. Moreover, given 
$g \in G$, we have $\tau \circ g^{*} = g^{*} \circ \tau$ 
on $H^{2}(X, \Bbb Z)$, because $\tau \circ g^{*} = g^{*} \circ \tau 
(= \sigma)$ on $S$ and $\tau \circ g^{*} = g^{*} \circ \tau 
(= g^{*})$ on $T$. Hence $\sigma \in O(S)^{+}$. \qed
\enddemo
\definition{(1.4)} Let us define the right object $\Cal N$ for $\Cal N'$. 
Unfortunately, $\Cal N$ is larger than the set $(\Cal N')^{G}$, in general.
Let $[b] \in \Cal N'$ and define a reduced divisor $B$ by 
$B := (\sum_{g \in G}g^{*}(b))_{red}$ and denote by
$B = \coprod_{k=1}^{n(b)} B_{k}$ the decomposition of $B$ into the connected 
components, where and in what follows, we identify $[b]$ 
and the unique smooth rational curve $b$ which represents $[b]$. 
Then, $G$ acts on the set $\{B_{k}\}_{k=1}^{n(b)}$ transitively, 
and the value $(B_{k}.B_{k})$ is then independent of $k$. 
Moreover, since $B_{k}$ is connected and reduced, 
using the Reimann-Roch Theorem and the exact sequence 
$0 \rightarrow \Cal O_{X}(-B_{k}) \rightarrow \Cal O_{X} \rightarrow 
\Cal O_{B_{k}} \rightarrow 0$, we see that  
$1 \geq 1 - p_{a}(B_{k}) = -(B_{k}.B_{k})/2$. This implies $(B_{k}.B_{k}) 
= -2$ if $(B_{k}.B_{k}) < 0$. Set 
$$\Cal N := \{b \in \Cal N' \vert (B_{k}.B_{k}) = -2\},$$ 
and define for $b \in \Cal N$ and $k \in \{1, 2, ..., n(b)\}$ the reflection 
$r_{B_{k}}$ on $H^{2}(X, \Bbb Z)$ 
by $r_{B_{k}}(x) = x + (x.B_{k})B_{k}$. It is easily seen that 
$r_{B_{k}}$ are Hodge isometries and satisfy $r_{B_{k}}(\Cal C') = \Cal C'$,
$r_{B_{k}} \circ r_{B_{l}} = r_{B_{l}} \circ r_{B_{k}}$ and 
$r_{B_{k}}^{2} = id$. 
\enddefinition 
By using the last two formulae in (1.4), we readily obtain 
\proclaim{Lemma (1.5)}
\roster
\item 
$(\prod_{k=1}^{n(b)}r_{B_{k}})(x) = x + \sum_{k=1}^{n(b)}(x.B_{k})B_{k}$ 
for $x \in H^{2}(X, \Bbb Z)$.
\item
$(\prod_{k = 1}^{n(b)} r_{B_{k}})^{2} = id$. \qed
\endroster
\endproclaim 
Now $G$-equivariant reflection group is defined as follows:
\definition{Definition (1.6)} Set 
$R_{b} := \prod_{k=1}^{n(b)} r_{B_{k}}$ for $b \in \Cal N$ and 
define a subgroup $\Gamma$ of $O(H^{2}(X, \Bbb Z))$ by  
$\Gamma := \,\langle R_{b}\, \vert\, b \in \Cal N \rangle$. 
We call this group $\Gamma$ the $G$-{\it equivariant reflection group 
of the pair} $(X, G)$. 
\enddefinition 
\proclaim{Lemma (1.7)} 
\roster
\item
$\Gamma \subset O(S)^{+}$, or more precisely, 
for each $\sigma \in \Gamma$, the restriction $\sigma \vert S$ lies in 
$O(S)^{+}$ and satisfies $\sigma = id$ if $\sigma \vert S = id$.
\item 
The cone $\Cal A$ is a fundamental domain for the action $\Gamma$ on 
$\Cal C$, that is, $\Cal A$ satisfies that 
$\sigma(\Cal A^{\circ}) \cap \Cal A^{\circ} \not= \emptyset$ if and only if 
$\sigma = id$, and that $\Gamma \cdot \Cal A = \Cal C$.
\endroster
\endproclaim 
\demo{Proof} Take $x \in S'$ and $g \in G$. Using (1.5), we calculate 
$$\align 
g^{*} \circ R_{b}(x) 
&= g^{*}(x) + \sum_{k=1}^{n(b)}(x.B_{k})g^{*}(B_{k})\\ 
&= g^{*}(x) + \sum_{k=1}^{n(b)}(g^{*}(x).g^{*}(B_{k}))g^{*}(B_{k})\\ 
&= g^{*}(x) + \sum_{k=1}^{n(b)}(g^{*}(x).B_{k})B_{k}\\ 
&= R_{b} \circ g^{*} (x).
\endalign$$ 
Hence $g^{*} \circ R_{b} = R_{b} \circ g^{*}$ 
and therefore, $R_{b}(S) = S$. Since $R_{b}(\Cal C') = \Cal C'$, 
this also gives $R_{b}(\Cal C) = \Cal C$. Moreover, 
$R_{b}$ is a Hodge isometry, 
because $R_{b}([\omega_{X}]) = [\omega_{X}]$. Therefore $R_{b} \in O(S)^{+}$. 
Assume that $\sigma \vert S =id$ 
for some $\sigma \in \Gamma$. Then, $\sigma(\Cal A^{\circ}) \cap 
\Cal A^{\circ} \not= 
\emptyset$, and in particular, $\sigma((\Cal A')^{\circ}) \cap 
(\Cal A')^{\circ} \not= \emptyset$. 
This implies $\sigma = id$, because $\Cal A'$ is the fundamental domain 
for the action $\langle r_{b}\, \vert\, b \in \Cal N' \rangle$ on $\Cal C'$ 
(see for example [BPV, Chap.VIII, Proposition (3.9)]). 
It remains to check the equality $\Gamma \cdot \Cal A = \Cal C$. 
Recall that $G$ acts transitively on the set 
$\{B_{k}\}_{k=1}^{n(b)}$ and satisfies $(x.g^{*}(B_{k})) = 
(g^{*}(x).g^{*}(B_{k})) = (x.B_{k})$ for $x \in S_{\Bbb R}$ 
if $g \in G$. Therefore, $(x.B_{k}) = (x.B_{l})$ 
for $x \in S_{\Bbb R}$ and for $b \in \Cal N$, and we get 
$R_{b}(x) = x + (x.B_{1})\sum_{k=1}^{n(b)}B_{k}$. 
This formula shows that $R_{b}\vert S_{\Bbb R}$ is nothing but  
a reflection with respect to the hyperplane defined by $(*.B_{1}) = 0$. 
Recall that by (1.3)(1) the subgroup (of index two) of the orthogonal group 
$O(S_{\Bbb R})$ 
preserving $\Cal C^{\circ}$ makes 
the quotient space $\Cal C^{\circ}/\Bbb R_{>0}$ a Lobachevskii 
space. Then, we can apply the general theory 
on the discrete reflection group on Lobachevskii space [Vi] 
to see that the space
$$\tilde{\Cal A} := 
\{x \in \Cal C\, \vert\, (x.B_{1}) \geq 0 \, 
\text{for all}\, b \in \Cal N \}$$ 
is a fundamental domain for the action $\Gamma$ on $\Cal C$. 
Therefore, it is sufficient 
to check $\tilde{\Cal A} = \Cal A$. 
It is clear that $\Cal A \subset \tilde{\Cal A}$. 
Let $x \in \tilde{\Cal A}$ and take $b \in \Cal N'$. 
Then, $(x.B_{1}) = c(x.b)$, where 
$c$ is the number of the irreducible components of $B_{1}$, 
because $x$ is $G$-invariant. 
If $b \not\in \Cal N$, then $(B_{1}.B_{1}) \geq 0$, and therefore,   
$(x.B_{1}) \geq 0$ by the Hodge index Theorem. If $b \in \Cal N$, then 
$(x.B_{1}) \geq 0$ by the definition of $\tilde{\Cal A}$. 
Hence, $(x.b) \geq 0$ for all $b \in \Cal N'$. 
This gives $\tilde{\Cal A} \subset \Cal A$. \qed
\enddemo 
Now we can formulate an equivariant Torelli Theorem as follows. 
This is a reformulation of the abstract version of the Torelli Theorem 
for K3 surface [SPP] in an equivariant setting and 
is also regarded as a sort of generalisation of the Torelli Theorem 
for Enriques surface due to Horikawa and Yukihiko Namikawa ([Ho], [Nm]).
\proclaim{Theorem (1.8)} $\Gamma$ is a normal subgroup of $O(S)^{+}$ and 
fits in with the semi-direct decomposition $O(S)^{+} = 
\Gamma \rtimes P(S)^{+}$. 
\endproclaim
\demo{Proof} Let $\sigma$ be an element of $O(S)^{+}$ 
and take $y \in \Cal A^{\circ}$. 
Applying (1.7)(2) for 
$\sigma (y)$, we find an element $r \in \Gamma$ such that 
$r^{-1} \circ \sigma(y) \in \Cal A^{\circ}$. Note that 
$r^{-1} \circ \sigma \in O(S)^{+}$ by (1.7)(1). 
Then, there exists a Hodge isometry 
$\rho \in O(H^{2}(X, \Bbb Z))$ such that $\rho \vert S = r^{-1} \circ \sigma$
and that $\rho \circ g^{*} = g^{*} \circ \rho$ for all $g \in G$. 
In addition, this $\rho$ is also effective, 
because $\rho(y) = r^{-1} \circ \sigma(y) \in \rho ((\Cal A')^{\circ}) 
\cap (\Cal A')^{\circ}$. Hence, by the Torelli Theorem for K3 
surfaces [PSS], [BPV, Chap.IIIV], there exists $f \in \text{Aut}(X)$ 
such that $f^{*} = \rho$. Moreover, $f \circ g = g \circ f$ 
for all $g \in G$ again by the Torelli Theorem, because 
$f^{*} \circ g^{*} = g^{*} \circ f^{*}$. Hence $f \in Q$ 
and $f^{*}\vert S \in P(S)^{+}$ by (1.3)(3). 
Since $f^{*}\vert S = r^{-1} \circ \sigma$, we get  
$O(S)^{+} = \Gamma \cdot P(S)^{+}$. 
Assume that $r \circ \tau = r' \circ \tau'$ for some
$r, r' \in \Gamma$ and $\tau, \tau' \in P(S)^{+}$. Since 
$\tau \circ (\tau')^{-1} (\Cal A) = \Cal A$ (1.3)(3), we have 
$r^{-1} \circ r'(\Cal A) = \Cal A$. Therefore, by (1.7)(2), 
we obtain $r^{-1} \circ r' = id$. This shows the uniqueness of 
the factorisation of elements of $O(S)^{+}$. 
It remains to show that $\Gamma$ is a normal subgroup of $O(S)^{+}$.  
For this, it is now enough to check that for each $b \in \Cal N$ and 
$\sigma \in P(S)^{+}$ there exists an element $b' \in \Cal N$ such that 
$\sigma^{-1} \circ R_{b} \circ \sigma = R_{b'}$. 
Let us choose $f \in Q$ such that $\sigma = f^{*} \vert S$. 
Then, we calculate 
$$\align
\sigma^{-1} \circ R_{b} \circ \sigma (x) 
&= (f^{-1})^{*} \circ R_{b} \circ f^{*} (x)\\ 
&= x + \sum_{i=1}^{n(b)} (f^{*}(x).B_{i})(f^{-1})^{*}(B_{i})\\ 
&= x + \sum_{i=1}^{n(b)} (x.(f^{-1})^{*}(B_{i}))(f^{-1})^{*}(B_{i}).
\endalign$$ 
In addition, we have $(f^{-1})^{*}(b) \in \Cal N$, because $b \in \Cal N$ 
and $f^{*} \circ g^{*} = g^{*} \circ f^{*}$.  
Therefore, we may take this $(f^{-1})^{*}(b)$ as $b'$.  
\qed
\enddemo 
Theorem (1.8) provides some more geometrical consequences. The first one is 
a generalisation of the main result of H. Sterk [St] (see also 
[Kaw5, Section 2]) to an equivariant setting:
\proclaim{Corollary (1.9)} 
There exists a finite rational polyhedral fundamental domain 
$\Delta$ for the action $P(S)^{+}$ on $\Cal A$.
\endproclaim 
\demo{Proof} Note that $O(S)^{+}$ is of finite index in 
the arithemetic group $O(S)$ of the self-dual 
homogeneous cone $\Cal C$ by (1.3)(4). Then, by [AMRT, Chap.II, 
Pages 116-117], there exists a 
finite rational polyhedral fundamental domain $\Delta$ for the action 
$O(S)^{+}$ on $\Cal C$. Translating $\Delta$ by 
an appropriate element of $\Gamma$ if necessary, we take such $\Delta$ as 
$(\Delta)^{\circ} \cap \Cal A^{\circ} \not= \emptyset$.  
This $\Delta$ satisfies $\Delta \subset \Cal A$. (Indeed, otherwise, there 
would be an element $b \in \Cal N$ such that the wall 
$H_{b} =\{x \in S_{\Bbb R} \vert (x.b)= 0\}$ of $\Cal A$ 
satisfies $(\Delta)^{\circ} \cap H_{b} \not= \emptyset$. 
However, since $R_{b}(y) = y$  for 
$y \in (\Delta)^{\circ} \cap H_{b}$, we would then have 
$R_{b}((\Delta)^{\circ}) \cap (\Delta)^{\circ} \not= \emptyset$, 
a contradiction.) Now, combining $P(S)^{+} = O(S)^{+}/\Gamma$ (1.8), 
$P(S)^{+}\cdot \Cal A = \Cal A$ (1.3)(3) and the fact that $\Cal A$ 
is a fundamental domain for the action $\Gamma$ on $\Cal C$ (1.7)(2), 
we conclude that this $\Delta$ gives a desired fundamental domain. 
\qed
\enddemo 
\proclaim{Corollary (1.10)} Let $Z$ be a normal K3 
surface and $G_{Z}$ 
a finite automorphism group 
of $Z$. Then $Z$ admits only finitely many $G_{Z}$-stable contractions 
up to $G_{Z}$-equivariant isomorphism. 
\endproclaim
\demo{Proof} First consider the case where $Z$ is smooth. 
Let $\Delta$ be the fundamental domain found in (1.9) for $(Z, G_{Z})$ and 
decompose $\Delta$ into its locally closed strata,   
$\Delta = \Delta_{1} \sqcup \Delta_{2} \sqcup \cdots \sqcup \Delta_{n}$. 
Then, any two $\Bbb Q$-rational points $H_{1}$ and $H_{2}$ 
in the same strata $\Delta_{i}$ give $G_{Z}$-equivariantly 
isomorphic $G_{Z}$-stable contractions, because, as there exist positive 
integers $m_{1}$, $m_{2}$ and $m_{3}$ such that $m_{2}H_{2} - m_{1}H_{1}, 
m_{3}H_{1} - m_{2}H_{2} \in \Delta_{i}$, by the semi-ampleness of rational 
nef divisors on a K3 surface [SD], the contractions given by 
$H_{1}$ and $H_{2}$ factor through each other and hence
$G_{Z}$-equivariantly isomorphic. 
Let $\Phi : Z \rightarrow W$ be a $G_{Z}$-stable contraction 
and choose 
a $G_{Z}$-invariant line bundle $H$ such that $\Phi = \Phi_{H}$. 
Then, by (1.9) and (1.3)(3), there exist an integer $i \in \{1, 2, ..., n\}$ 
and $f \in Q$ (for $Z$) such that 
$f^{*}([H]) \in \Delta_{i}$. Therefore, the result follows for smooth $Z$. 
\par
Next consider the case where $Z$ is singular. Let $\nu : Y \rightarrow Z$ be 
the minimal resolution, 
$G_{Y}$ the unique equivariant lift of the action $G_{Z}$ on $Y$ and 
$E$ the exceptional divisor of $\nu$. 
Then, $E$ is a disjoint sum of reduced divisors of Dynkin type 
and is stable under $G_{Y}$. 
Moreover, $Y$ admits only finitely many 
$G_{Y}$-stable contractions up 
to $G_{Y}$-equivariant isomorphism by the previous argument.
Let us denote the representatives of the $G_{Y}$-stable contractions 
of $Y$ by 
$\Phi_{i} : Y \rightarrow W_{i}$ ($i = 1, 2, ..., I$). 
Each of these $\Phi_{i}$ is either 
an elliptic fibration or a birational contraction. In particular,
for each $\Phi_{i}$, there are only finitely many effective reduced divisors 
$E_{ij}$ ($j = 1, 2, ..., J_{i}$) on $Y$ satisfying that 
$E_{ij}$ is a disjoint sum of reduced divisors of Dynkin type, 
stable under $G_{Y}$ and 
$\text{dim}\,\Phi_{i}(E_{ij}) = 0$. 
Let $\nu_{ij} : Y \rightarrow Z_{ij}$ be the contraction of $E_{ij}$ and 
$\varphi_{ij} : Z_{ij} \rightarrow W_{i}$ the induced contraction. 
Then, $G_{Y}$ descends 
equivariantly to the action on 
$Z_{ij}$ and makes $\varphi_{ij}$ $G_{Y}$-stable.
Let $\varphi : Z \rightarrow W$ be a $G_{Z}$-stable contraction of $Z$ and set 
$\Phi := \varphi \circ \nu : Y \rightarrow W$. 
Then, this $\Phi$ is $G_{Y}$-stable. Therefore, there exist $i$ and 
isomorphisms $F : Y \simeq Y$ and 
$f : W \simeq W_{i}$ 
such that $\Phi : Y \rightarrow W$ and 
$\Phi_{i} : Y \rightarrow W_{i}$ are 
$G_{Y}$-equivariantly isomorphic by $F$ and $f$. Since $F(E)$ satisfies all 
the defining 
properties of $E_{ij}$, there also exists $j$ such that $F(E) = E_{ij}$. 
Therefore, $F$ descends to give an isomorphism $F' : Z \rightarrow Z_{ij}$.  
This $F'$ together with $f$ gives a $G_{Z}$-equivariant isomorphism 
(with respect to the representation $G_{Z} \rightarrow \text{Aut}(Z_{ij})$ 
through $G_{Y}$) between $\phi : Z \rightarrow W$ and 
$\varphi_{ij} : Z_{ij} \rightarrow W_{i}$. Now we are done. \qed
\enddemo 
Next, as an application of (1.7)(2), we show the following generalisation of 
the result of S. Kondo [KS, Lemma(2.1)]. 
\proclaim{Corollary (1.11)} Let $(X, G)$ be a pair of a K3 surface and its 
finite automorphism group and, as before, $S$ the $G$-invariant part of the 
N\'eron-Severi group of $X$.
\roster
\item
Assume that $S$ represents $0$. Then $X$ admits a $G$-stable elliptic 
fibration. In particular, if the rank of $S$ is greater than or equal to $5$, 
then $X$ admits a $G$-stable elliptic fibration.
\item
Assume that $S$ contains the even unimodular hyperbolic lattice 
$U$ of rank $2$. Then $X$ admits a $G$-stable Jacobian fibration, 
that is, a $G$-stable elliptic fibration having at least one $G$-stable 
global sections. In particular, if $\text{rank}\, S \geq 3 + l(S)$, 
where $l(S)$ is the minimal number of generaters of the finite abelian group 
$S^{*}/S$, then $X$ admits a $G$-stable Jacobian fibration. 
\endroster
\endproclaim
\demo{Proof of (1)} By the assumption, there exists a primitive point 
$x \in \partial \Cal C \cap S$. By (1.7)(2), translating $x$ by 
an appropriate element of $\Gamma$, we obtain a primitive point 
$y \in \Cal A$ such that $(y.y) = 0$. 
This $y$ gives a $G$-stable elliptic fibration on 
$X$. The last assertion now follows from the famous arithmetical fact due to 
[Se, Page 43] that every indefinite rational quadratic form of 
$n$-variables represents 0 provided that $n\geq 5$.
\qed 
\enddemo
\demo{Proof of (2)} Choose an integral basis $\langle e, c \rangle$ 
of $U$ such that $(e.e) = 0$, $(e.c) = 1$ and $(c.c) = -2$. 
Then as in (1), by translating $e$ by an appropriate element of $\Gamma$, 
we may assume that $e$ is the class of a smooth fiber $E$ of an 
$G$-stable elliptic fibration $\varphi : X \rightarrow \Bbb P^{1}$.  
Let us also choose an divisor $C$ such that $c = [C]$. Then, 
by the Riemann-Roch Theorem and the Serre duality, we obtain 
$$h^{0}(\Cal O_{X}(C)) + h^{0}(\Cal O_{X}(-C)) \, \geq 
\chi(\Cal O_{X}(C)) = (c^{2}/2) + 2 = 1.$$ 
Since $(C.E) = 1$, we have $h^{0}(\Cal O_{X}(-C)) = 0$.  
Therefore $h^{0}(\Cal O_{X}(C)) > 0$. Let $\vert C \vert = F + \vert M \vert$ 
be the decomposition of $\vert C \vert$ into the fixed component and the 
movable part. Since $\vert M \vert$ is free ([SD, Corollary 3.2]), 
by the Bertini Theorem, we may assume that $M$ itself is a smooth divisor. 
Note also that the divisor $F$ is $G$-stable by $g^{*}\vert C \vert 
= \vert C \vert$ for all $g \in G$ and by the uniqueness of the fixed 
part. In addition, we have either $(F.E) = 0$ and $(M.E) = 1$ or 
$(F.E) = 1$ and $(M.E) = 0$, because 
$1 = (C.E) = (F.E) + (M.E)$ and $E$ is nef. 
Assume that the first case happens. 
Then, $M$ must be a global section of $\varphi$ 
by the smoothness of $M$ and by $(E.M) = 1$, 
and would then satisfy $(M^{2}) = -2$. 
However, this contradicts the fact that $\vert M \vert$ is the movable 
part. 
Therefore, $(M.E) = 0$ and $(F.E) = 1$. 
Write the irreducible 
decomposition of $F$ as $F = C_{0} + \sum m_{i}C_{i}$, where 
$(C_{0}.E) = 1$ and 
$(C_{i}.E) = 0$ for $i \geq 1$. Then $C_{0}$ is a section and 
all other $C_{i}$ ($i \geq 1$) are vertical with respect to $\varphi$. 
Thus, $C_{0}$ is also $G$-stable. Therefore, this $C_{0}$ 
gives a desired section. The last statement now follows 
from the splitting Theorem due to Nikulin [Ni2, Corollary 1.13.5]. 
\qed 
\enddemo
\head
{\S 2. Calabi-Yau threefolds with infinite fundamental group}
\endhead
\par
\vskip 4pt
In this section we study Calabi-Yau threefolds with infinite 
fundamental group. As is already remarked in Introduction, such threefolds 
are smooth. Therefore, 
we may speak of their minimal splitting coverings introduced by Beauville, 
which, in threefold case, is summarised as follows: 
\definition{Summary (2.1)([Be2, Section 3])} 
Let $X$ be a smooth threefold with infinite 
fundamental 
group such that $c_{1}(X) = 0$ in $H^{2}(X, \Bbb R)$. Then, by the 
Bogomolov decomposition Theorem, such an $X$ admits an \'etale Galois 
covering either by an abelian threefold or by the product of a K3 surface 
and an elliptic curve. We also call $X$ {\it of Type A} in the former case 
and {\it of Type K} in the latter case.
Among many candidates of such coverings for a given 
$X$, there always exists the smallest one which is known to be unique 
for each $X$ up to isomorphism as a covering space and is obtained by posing 
one additional condition on the Galois group $G$ that $G$ contains no 
non-zero translations in the former case and that $G$ contains no elements 
of the form $(id, \text{non-zero translation})$ in the later case. 
According to Beauville, 
we call this smallest covering {\it the minimal splitting covering} of $X$. 
\qed
\enddefinition

We study Calabi-Yau threefolds of Type A in the subsection (2A)
and those of Type K in (2K) through their minimal 
aplitting coverings. 
\head
{\S 2A. Calabi-Yau threefolds of Type A}
\endhead
\par
\vskip 4pt
The aim of this subsection is to show Theorem (0.1) in Introduction. 

\demo{Proof of (0.1)(I)} 
\par
\vskip 4pt 
First of all, we introduce the following:  
\definition{Definition (2.2)} A finite group $G$ is called a Calabi-Yau group 
of Type A (resp. a pre-Calabi-Yau group of Type A), which, 
throughout this subsection, is abbreviated by a C.Y. group (resp. by a 
pre-C.Y. 
group), if there exist an abelian 
threefold $A$ and a faithful representation $G \hookrightarrow \text{Aut}(A)$ 
which satisfy the following conditions (1) - (4) (resp. (1) - (3)):
\roster
\item 
$G$ contains no non-zero translations;
\item 
$g^{*}\omega_{A} = \omega_{A}$ for all $g \in G$;
\item 
$A^{[G]} = \emptyset$;
\item 
$H^{0}(A, \Omega_{A}^{1})^{G} = \{0\}$. 
\endroster
In each case we call $A$ a target abelian threefold. \qed
\enddefinition

Then the proof of (I) is equivalent to classifying C.Y. groups 
together with their actions on target abelian threefolds. The 
following inductive nature of pre-C.Y. groups turns out to be useful:

\proclaim{Lemma (2.3)}  If $G$ is a pre-C.Y. group, then so are all the 
subgroups of $G$. In other words, if $G$ contains a non pre-C.Y. group, 
then $G$ is not a pre-C.Y. group, either. \qed \endproclaim

\proclaim{Lemma (2.4)} Let $G$ be a pre-C.Y. group, $A$ its target abelian 
threefold 
and $\rho : G \rightarrow \text{GL}(H^{0}(A, \Omega_{A}^{1}))$ 
the natural representation. Then:
\roster
\item 
$\rho$ is injective.
\item 
$\text{Im}(\rho) \subset \text{SL}(H^{0}(A, \Omega_{A}^{1}))$.
\item 
Let $g$ be an element of $G$ and set $n = \text{ord}(g)$. Then, $n \in 
\{1, 2, 3, 4, 6\}$. Moreover, there exists a basis of 
$H^{0}(A, \Omega_{A}^{1})$ under which 
$g^{*} \vert H^{0}(A, \Omega_{A}^{1}) = 
\text{diag}(1, \zeta_{n}^{k}, \zeta_{n}^{-k})$ 
for some $k$ such that $(n, k) = 1$. 
\endroster 
\endproclaim 

\demo{Proof} The assertions (1) and (2) are clear. 
Let us show the assertion (3). 
Choose a basis of $H^{0}(A, \Omega_{A}^{1})$ under which 
the matrix representaion of $g^{*} \vert H^{0}(A, \Omega_{A}^{1})$ is 
diagonalised and write $g^{*} = \text{diag}(a, b, c)$. 
Then, there exist global coordinates $(x, y, z)$ of $A$ 
such that the (co-)action of $g$ on $A$ is of 
the form  $g(x, y, z) = (ax + p, by + q, cz + r)$.
Suppose $a \not= 1$, $b \not= 1$, and $c \not= 1$. Then, 
the point $P = (p/(1-a), q/(1-b), r/(1-c)) \in A$ is a fixed 
point of $g$, a contradiction. Therefore we may assume $a = 1$, 
$b = \zeta_{n}$ and $c = \zeta_{-n}$ (by replacing $g$ by 
an appropriate generater of $\langle g \rangle$ and reordering 
the basis if necessary). Recall that $H^{1}(A, \Bbb C) 
= H^{0}(A, \Omega_{A}^{1}) \oplus 
\overline{H^{0}(A, \Omega_{A}^{1})}$. Then, 
$g^{*} \vert H^{1}(A, \Bbb C) = \text{diag}(1, \zeta_{n}, \zeta_{n}^{-1}, 1, 
\zeta_{n}^{-1}, \zeta_{n})$. 
Hence, $\varphi(n) \leq (6-2)/2 = 2$, and therefore, 
$n \in \{1, 2, 3, 4, 6\}$, 
because $g^{*} \vert H^{1}(A, \Bbb C)$ is the scalar 
extension of $g^{*} \vert H^{1}(A, \Bbb Z)$. 
\qed \enddemo

First, we determine commutative pre-C.Y. groups.

\proclaim{Lemma (2.5)} Let $G$ be a commutative pre-C.Y. group and $A$ a 
target abelian threefold. Then $G$ is isomorphic to either $C_{n}$, 
where $1 \leq n \leq 6$ and $n \not= 5$, or $C_{2}^{\oplus 2}$. 
In particular, there 
exist no commutative pre-C.Y. groups of order $\geq 7$. 
Moreover, if $G$ is a commutative C.Y. group, then $G$ is isomorphic to 
$C_{2}^{\oplus 2}$ and the action of $G$ on $H^{0}(A, \Omega_{A}^{1})$ is 
same as in (0.1) (I) (1). 
\endproclaim

\demo{Proof} Set 
$G = \langle g_{1} \rangle \oplus 
\cdots \oplus \langle g_{r} \rangle \simeq C_{n_{1}} \oplus C_{n_{2}} \oplus 
\cdots \oplus C_{n_{r}}$, 
where $r \geq 0$ and $2 \leq n_{1} \vert n_{2} \vert \cdots \vert n_{r}$. 
If $r \leq 1$, then we get the result by (2.4)(3). Assume that 
$r \geq 2$. Let $i, j$ be two integers such that $1 \leq i < j \leq r$. 
Then, using $g_{i}g_{j} = g_{j}g_{i}$ and (2.4)(3), and 
replacing $g_{i}$ and $g_{j}$ by 
other generaters of $\langle g_{i} \rangle$ and $\langle g_{j} \rangle$ if 
necessary, we may find a basis of $H^{0}(A, \Omega_{A}^{1})$ under which 
$g_{i}^{*} \vert H^{0}(A, \Omega_{A}^{1})$ and $g_{j}^{*} \vert H^{0}(A, 
\Omega_{A}^{1})$ are simulteneously diagonalised as either 
one of the following forms:

(1) $g_{i}^{*} = \text{diag}(1, \zeta_{n_{i}}, \zeta_{n_{i}}^{-1})$ and 
$g_{j}^{*} = \text{diag}(1, \zeta_{n_{j}}, \zeta_{n_{j}}^{-1})$ or 

(2) $g_{i}^{*} = \text{diag}(1, \zeta_{n_{i}}, \zeta_{n_{i}}^{-1})$ and 
$g_{j}^{*} = \text{diag}(\zeta_{n_{j}}^{-1}, 1, \zeta_{n_{j}})$. 

In the former case, we have $g_{i}^{*} = (g_{j}^{*})^{n_{j}/n_{i}}$, 
whence by (2.4)(1), $g_{i} = (g_{j})^{n_{j}/n_{i}}$, a contradiction. 
In the latter case, 
we calculate 
\par
$g_{i}^{*} g_{j}^{*} = \text{diag}(\zeta_{n_{j}}^{-1}, 
\zeta_{n_{i}}, \zeta_{n_{i}}^{-1}\zeta_{n_{j}})$ and 
$(g_{i}^{-1})^{*}g_{i}^{*} = 
\text{diag}(\zeta_{n_{j}}^{-1}, \zeta_{n_{i}}^{-1}, 
\zeta_{n_{i}}\zeta_{n_{j}})$. 
\par
Therefore, by(2.4)(3), we get 
$\zeta_{n_{i}}^{-1}\zeta_{n_{j}} = \zeta_{n_{i}}\zeta_{n_{j}} = 1$. 
This implies $\zeta_{n_{i}} = \zeta_{n_{j}} = -1$ and $n_{i} = 2$ for all 
$i = 1, 2, \cdots, r$ again by (2.4)(3). In particular,
\par 
$g_{1}^{*} = \text{diag}(1, -1, -1)$ and 
$g_{2}^{*} = \text{diag}(-1, 1, -1)$.
\par 
Assume that $r \geq 3$. 
Then, $g_{3}^{*}$ must be of the form $\text{diag}(-1, -1, 1)$. 
However, then, $g_{1}^{*}g_{2}^{*} = g_{3}^{*}$, 
whence $g_{1}g_{2} = g_{3}$, a contradiction. 
Therefore, if $r \geq 2$, then $G \simeq C_{2}^{\oplus 2}$ and 
there exists a basis of $H^{0}(A, \Omega_{A}^{1})$ under which $g_{1}^{*} = 
\text{diag}(1, -1, -1)$ and $g_{2}^{*} = \text{diag}(-1, 1, -1)$. 
The remaining assertions follow from this description and (2.4)(3). 
\qed 
\enddemo

Let us examine non-comutative pre-C. Y. groups. First we estimate their 
orders. The following two Theorems are extremally useful:

\proclaim{Theorem (2.6) (Wielandt, eg.[KT, Chap.2, Theorem (2.2)])} 
Let $G$ be an arbitray finite group, $p$ a prime number and $h$ a positive 
integer such that $p^{h} \vert \vert G \vert$. Then there exists a subgroup 
$H$ of $G$ such that $\vert H \vert = p^{h}$. \qed \endproclaim

\proclaim{Theorem (2.7) (Burnside-Hall, eg. [Su, Page 90, Corollary 2])} 
Let $K$ be an arbitray $p$-group and $H$ a maximal, normal commutative 
subgroup of $G$. Set $\vert G \vert = p^{n}$ and $\vert H \vert = p^{h}$. 
Then $h(h+1) \geq 2n$. \qed \endproclaim 

\proclaim{Lemma (2.8)} Let $G$ be a pre-C.Y. group. Then, 
$\vert G \vert$ is either $2^{n}$ or $2^{n}\cdot 3$, where $n$ 
is an integer such that $0 \leq n \leq 3$, or more explicitly, 
$\vert G \vert \in \{1, 2, 3, 4, 6, 8, 12, 24\}$. 
\endproclaim 

\demo{Proof} By (2.6) and (2.4)(3), we have $\vert G \vert = 
2^{n}\cdot3^{m}$, where $m$ and $n$ are some non-negative integers. 
Assume for a contradiction that $m \geq 2$. 
Then, by(2.6), $G$ contains a subgroup $H$ of order $3^{2}$.
This $H$ must be a pre-C.Y. group by (2.3) and is also commutative. 
However, this contradicts (2.5). 
Therefore $m = 0$ or $1$. Assume to the contrary, that $n \geq 4$. 
Then, again by (2.6), $G$ contains 
a subgroup $H$ of order $2^{4}$. Let $K$ be a maximal normal commutative 
subgroup of $H$ and set $\vert K \vert = 2^{k}$. By (2.7), 
this $k$ satisfies $k(k+1) \geq 8$ and hence $k \geq 3$. 
However, this again contradicts (2.3) and (2.5). 
Therefore $n \leq 3$. \qed \enddemo

Combining this with the classification of non-commutative finite groups 
of small order (eg. [Bu, Chap.4, Pages 54-55 and Chap.5, Pages 83-89]), 
we get: 

\proclaim{Corollary (2.9)} 
Let $G$ be a pre-C.Y. group. Assume that $G$ is non-commutative 
and satisfies $\vert G \vert \leq 12$. Then $G$ is isomorphic to either one of
$D_{6} (\simeq S_{3})$, $D_{8}$, $Q_{8}$, $D_{12}$, $Q_{12}$ or $A_{4}$. 
\qed \endproclaim 

Next we show that among the candidates in (2.9), only $D_{8}$ 
can be realised as a C.Y. group of Type A. For proof, let us recall 
the following:

\proclaim{Proposition (2.10) (eg. [KT, Chap.8, Pages 273-275])}
Up to equivalence, the complex linear irreducible representaions of $D_{2n}$ 
($3 \leq n \in \Bbb Z$), $Q_{4n}$ ($1 \leq n \in \Bbb Z$) and $A_{4}$ are 
given as follows: 
\flushpar
$(D_{0})$.  $D_{2n} = \langle a, b \vert a^{n} = b^{2} = 1, bab = a^{-1} 
\rangle$ such that $n \equiv 0 (\text{mod}\, 2)$: 
\roster
\item 
$\rho_{1, 0} : a \mapsto 1, b \mapsto 1$; $\rho_{1, 1} : a \mapsto 1, b 
\mapsto -1$; $\rho_{1, 2} : a \mapsto -1, b \mapsto 1$; $\rho_{1, 3} : a 
\mapsto -1, b \mapsto -1$;
\item 
$\rho_{2, k} : a \mapsto \left(\matrix \zeta_{n}^{k} & 0 \\ 0 & 
\zeta_{n}^{-k}\endmatrix \right), b \mapsto \left(\matrix 0 & 1 \\ 1 
& 0 \endmatrix \right)$, where $k$ is an integer such that $1 \leq k \leq n/2 
- 1$.
\endroster
\flushpar
$(D_{1})$.  $D_{2n} = \langle a, b \vert a^{n} = b^{2} = 1, bab = a^{-1} 
\rangle$ such that $n \equiv 1 (\text{mod}\, 2)$: 
\roster
\item 
$\rho_{1, 0} : a \mapsto 1, b \mapsto 1$; $\rho_{1, 1} : a \mapsto 1, b 
\mapsto -1$; 
\item 
$\rho_{2, k} : a \mapsto \left(\matrix \zeta_{n}^{k} & 0 \\ 0 & 
\zeta_{n}^{-k}\endmatrix \right), b \mapsto \left(\matrix 0 & 1 \\ 1 
& 0 \endmatrix \right)$, where $k$ is an integer such that $1 \leq k \leq 
(n - 1)/2$.
\endroster
\flushpar
$(Q_{0})$. $Q_{4n} = \langle a, b \vert a^{2n} = 1, a^{n} = b^{2}, b^{-1}ab = 
a^{-1} \rangle$ such that $n \equiv 0 (\text{mod}\, 2)$:
\roster
\item 
$\rho_{1, 0} : a \mapsto 1, b \mapsto 1$; $\rho_{1, 1} : a \mapsto 1, b 
\mapsto -1$; $\rho_{1, 2} : a \mapsto -1, b \mapsto 1$; $\rho_{1, 3} : a 
\mapsto -1, b \mapsto -1$;
\item 
$\rho_{2, l} : a \mapsto \left(\matrix \zeta_{2n}^{l} & 0 \\ 0 & 
\zeta_{2n}^{-l}\endmatrix \right), b \mapsto \left(\matrix 0 & \zeta_{4} \\ 
\zeta_{4} & 0 \endmatrix \right)^{l}$, where $l$ is an integer such that $1 
\leq k \leq n - 1$.
\endroster
\flushpar
$(Q_{1})$. $Q_{4n} = \langle a, b \vert a^{2n} = 1, a^{n} = b^{2}, b^{-1}ab = 
a^{-1} \rangle$ such that $n \equiv 1 (\text{mod}\, 2)$:
\roster
\item 
$\rho_{1, 0} : a \mapsto 1, b \mapsto 1$; $\rho_{1, 1} : a \mapsto 1, b 
\mapsto -1$; $\rho_{1, 2} : a \mapsto -1, b \mapsto \zeta_{4}$; $\rho_{1, 3} 
: a \mapsto -1, b \mapsto -\zeta_{4}$;
\item 
$\rho_{2, l} : a \mapsto \left(\matrix \zeta_{2n}^{l} & 0 \\ 0 & 
\zeta_{2n}^{-l}\endmatrix \right), b \mapsto \left(\matrix 0 & \zeta_{4} \\ 
\zeta_{4} & 0 \endmatrix \right)^{l}$, where $l$ is an integer such that $1 
\leq k \leq n - 1$.
\endroster
\flushpar
$(A_{4})$. $A_{4}  = \langle a, b \rangle (\subset S_{4})$, 
where $a = (1 2 3)$ and $b = (1 2)(3 4)$: 
\roster
\item 
$\rho_{1, 0} : a \mapsto 1, b \mapsto 1$; $\rho_{1, 1} : a \mapsto \zeta_{3}, 
b \mapsto 1$; $\rho_{1, 2} : a \mapsto \zeta_{3}^{-1}, b \mapsto 1$; \item 
$\rho_{3} : a \mapsto \left(\matrix 0 & 1 & 0 \\ 0 & 0 & 1 \\ 1 & 0 & 0 
\endmatrix \right), b \mapsto \left(\matrix 1 & 0 & 0 \\ 0 & -1 & 0 
\\ 0 & 0 & -1 \endmatrix \right)$. \qed 
\endroster
\endproclaim

\proclaim{Lemma (2.11)} Let $G$ be a pre-C.Y. group and $A$ a target 
abelian threefold. Assume that 
$G$ is isomorphic to either $D_{8}$, $Q_{8}$ or $Q_{12}$. 
Then, the irreducible decomposition of the natural representation 
$\rho : G \rightarrow 
SL(H^{0}(A, \Omega_{A}^{1}))$ is of the form $\rho = \rho_{1, 1} \oplus 
\rho_{2,1}$ if $G \simeq D_{8}$ and $\rho = \rho_{1, 0} \oplus \rho_{2,1}$ if 
$G \simeq Q_{8}$ or $Q_{12}$, where we adopt the same notation as 
in (2.10). In particular, if $G$ is a C.Y. group, then 
$G \simeq D_{8}$ and the representation of $G$ on $H^{0}(A, \Omega_{A}^{1})$ 
is equivalent to the one given in (0.1) (I) (2). \endproclaim

\demo{Proof} Note that $\rho$ is not isomorphic to a direct sum of three 
1-dimensional representations, because $\rho$ is injective 
and $G$ is non-commutative. Then, the result follows from 
the list in (2.10) and (2.4)(1),(2). 
\qed 
\enddemo

The next two Lemmas are crucial and their proofs 
involve geometric consideration based on the non-cohomological 
condition $A^{[G]} = \emptyset$. 

\proclaim{Lemma (2.12)} Neither $D_{6} (\simeq S_{3})$ nor $D_{12}$ 
is a pre-C.Y. group. \endproclaim

\demo{Proof} The assertion for $D_{12}$ follows from the one for $D_{6}$ 
by (2.3), because $D_{6}$ can be embedded in $D_{12}$. 
Assume to the contrary, that $D_{6} = \langle a, b \vert a^{3} = b^{2} = 1, 
bab = a^{-1} \rangle$ 
is a pre-C.Y. group. Let $A$ be a target abelian threefold and 
$\rho : D_{6} \rightarrow SL(H^{0}(A, \Omega_{A}^{1}))$ the natural 
representation. Then, by the same argument as in (2.11), we obtain $\rho = 
\rho_{1, 1} \oplus \rho_{2,1}$. In other words, there exists a basis 
$\langle v_{1}, v_{2}, v_{3} \rangle$ of $H^{0}(A, \Omega_{A}^{1})$ under 
which $a^{*} = \left(\matrix 1 & 0 & 0 \\ 0 & \zeta_{3} & 0 \\ 0 & 0 & 
\zeta_{3}^{-1} \endmatrix \right)$ and $b^{*} = \left(\matrix -1 & 0 & 0 \\ 0 
& 0 & 1 \\ 0 & 1 & 0 \endmatrix \right)$. Let us fix an origin $0$ of $A$ and 
regard $A$ as a group variety. Set $\alpha := a(0)$ and $\beta := b(0)$. 
Then, we have $a = t_{\alpha}\circ a_{0}$ and $b = t_{\beta} \circ b_{0}$, 
where $a_{0}$, $b_{0}$ are the Lie part of $a$, $b$. 
Set $\tilde{E} := \text{Ker}(a_{0} - id_{A} : A \rightarrow A)$. This is 
a one-dimensional subgroup scheme of $A$. Let us take the identity component 
$E$ of $\tilde{E}$ and consider the quotient homomorphism  
$\pi : A \rightarrow S := A/E$. Notice that $E$ is an elliptic curve, 
$S$ is an abelian surface and the fibers of $\pi$ are of the form 
$E + s$ ($s \in A$). 
\proclaim{Claim (2.13)} $G$ descends to an automorphism group of $S$, that is, 
there exist automorphisms $\overline{a}$ and $\overline{b}$ 
of $S$ such that $\overline{a} \circ \pi = \pi \circ a$ and 
that $\overline{b} \circ \pi = \pi \circ b$. 
\endproclaim
\demo{Proof of Claim} Let $F$ be a fiber of $\pi$ and write $F = E + s$.  
Note that for $x \in A$, we have 
$a(x+s) = t_{\alpha}(a_{0}(x+s)) = t_{\alpha}(a_{0}(x)+ a_{0}(s)) 
= a_{0}(x)+ (a_{0}(s) + \alpha)$. 
Since $a_{0}(E) = E$, this formula implies $a(E+s) = E+ (a_{0}(s) + \alpha)$.  
Therefore, $a$ descends to 
the automorphism $\overline{a}$ of $S$ given $\pi(s) \mapsto 
\pi(a_{0}(s) + \alpha)$. Similarly, we have $b(x + s) = b_{0}(x) + 
(b_{0}(s) + \beta)$. Moreover, using $a_{0}b_{0} = b_{0}a_{0}^{-1}$, we 
calculate $a_{0}(b_{0}(e)) = b_{0}(a_{0}^{-1}(e)) = b_{0}(e)$ 
for $e \in E$. Therefore $b_{0}(E) \subset \tilde{E}$.
This implies $b_{0}(E) = E$, because $b_{0}(0) = 0 \in E$. 
Hence, $b(E+s) = E + (b_{0}(s) + \beta)$, and $b$ also descends 
to the automorphism $\overline{b}$ of $S$ defined 
by $\pi(s) \mapsto \pi(b_{0}(s) + \beta)$. \qed 
\enddemo
By construction, there exists a basis $\langle \overline{v}_{2}, 
\overline{v}_{3} \rangle$ of $H^{0}(S, \Omega_{S}^{1})$ such that 
$\pi^{*}(\overline{v}_{2}) 
= v_{2}$ and $\pi^{*}(\overline{v}_{3}) = v_{3}$. Using this basis, 
we have 
$\overline{a}^{*} = \left(\matrix \zeta_{3} & 0 \\ 0 & \zeta_{3}^{-1} 
\endmatrix \right)$ and $\overline{b}^{*} = \left(\matrix 0 & 1 \\ 1 & 0 
\endmatrix \right)$ on $H^{0}(S, \Omega_{S}^{1})$. 
This expression, in particular, shows that $S^{\overline{a}}$ 
consists of isolated points. Therefore, using the canonical graded isomorphism 
$$H^{*}(S, \Bbb C) = \oplus_{k=0}^{4} \wedge^{k}(H^{0}(S, \Omega_{S}^{1}) 
\oplus \overline{H^{0}(S, \Omega_{S}^{1})})$$ 
and applying topological Lefschetz fixed point formula, we readily obtain that 
$\vert S^{\overline{a}} \vert  = 9$.  
On the other hand, the equality $ab = ba^{-1}$ gives an equality 
$\overline{a}\overline{b} = \overline{b}\overline{a}^{-1}$. 
Therefore, $\overline{b}$ acts on the nine point 
set $S^{\overline{a}}$ 
and has a fixed point $\overline{s} \in S^{\overline{a}}$, 
because $\text{ord}(\overline{b}) = 2$. Put $F := \pi^{-1}(\overline{s})$. 
Then, $b(F) = F$ and $b^{*} \vert H^{0}(F, \Omega_{F}^{1}) = -1$ by 
the description of $a$ and $b$. Since $F$ is an elliptic curve, $F^{b}$ 
would then be non-empty. However this contradicts $A^{[G]} = \emptyset$. \qed 
\enddemo

\proclaim{Lemma (2.14)} The group $A_{4}$ is not a pre-C.Y. group. 
\endproclaim
\demo{Proof} Assume to the contrary that $A_{4} = \langle a, b \rangle$ is a 
pre-C.Y. group, where $a$ and $b$ denote the elements defined in (2.10). Let 
$A$ be a target abelian threefold and express $A$ as 
$A = \Bbb C^{3}/\Lambda$, where $\Lambda$ is a discrete sublattice of 
$\Bbb C^{3}$ of rank $6$. 
(In this proof, we 
regard $A$ as a three dimensional complex torus rather than an abelian 
variety.) Then, using the same argument as in (2.11), we readily find 
global coordinates $(z_{1}, z_{2}, z_{3})$ of $A$ under which 
the (co-)actions of $a$ and $b$ on $A$ 
are written as follows: 
$$a(z_{1}, z_{2}, z_{3}) = (z_{2}, z_{3}, z_{1}) + (\alpha_{1},\alpha_{2},
\alpha_{3}),\, b(z_{1}, z_{2}, z_{3}) = (z_{1}, -z_{2}, -z_{3}) + 
(\beta_{1},\beta_{2},\beta_{3}).$$ 
By this description, we obtain $a^{3}(z_{1}, z_{2}, z_{3}) 
= (z_{1}, z_{2}, z_{3}) + 
(\alpha,\alpha,\alpha)$, 
where we put $\alpha := \alpha_{1}+\alpha_{2}+\alpha_{3}$. 
Since $a^{3} = id$, we have $(\alpha,\alpha,\alpha) \in \Lambda$. Set 
$t := t_{(\alpha, \alpha, \alpha)}$. Then, on the one hand,  
$b^{-1}\circ t \circ b(z_{1}, 
z_{2}, z_{3}) = (z_{1}, z_{2}, z_{3}) + (\alpha,-\alpha,-\alpha)$ 
and, on the other hand, $b^{-1}\circ t \circ b = id$, because 
$(\alpha, \alpha, \alpha) \in \Lambda$. Therefore, 
$(\alpha,-\alpha,-\alpha) \in \Lambda$ and hence, 
$(2\alpha, 0, 0) = (\alpha,\alpha,\alpha) + (\alpha,-\alpha,-\alpha) 
\in \Lambda$. Consider the point $P := [(0, \alpha_{2} + \alpha_{3}, 
\alpha_{1} + \alpha_{2} + 2\alpha_{3})] \in A$. Then, 
$a^{2}(P) 
= (2\alpha_{1} + 2\alpha_{2} + 2\alpha_{3}, \alpha_{2} + 2 \alpha_{3}, 
\alpha_{1} + \alpha_{2} + 2\alpha_{3}) 
= (0, \alpha_{2} + \alpha_{3}, \alpha_{1} + \alpha_{2} + 2 \alpha_{3}) 
+ (2 \alpha, 0, 0) = P$, a contradiction. 
\qed 
\enddemo 
In order to complete the proof of (0.1) (I), it remains to show the following:

\proclaim{Lemma (2.15)} Let $G$ be a group of order $24$. Then $G$ is not 
a C.Y. group. \endproclaim 

\demo{Proof} We shall assume to the contrary, that there exists a C.Y. group 
$G$ of order $24$, and derive a contradiction. Our argument is based 
on the following classification of the groups of order $24$:

\proclaim{Proposition (2.16) (eg. [Bu, Chapter 9, Pages 171-174])} 
Let $G$ be an (arbitrary) group of order $24$, $H_{2}$ a $2$-Sylow subgroup 
of $G$ and $H_{3} = \langle c \rangle (\simeq C_{3})$ a $3$-Sylow subgroup of 
$G$. Then, $H_{2}$ is isomorphic to either $C_{8}$, $C_{2} \oplus C_{4}$, 
$C_{2}^{\oplus 3}$, $D_{8}$ or $Q_{8}$ and $G$ is isomorphic to one of the 
following 15 groups according to the isomorphism class of 
the $2$-Sylow subgroup $H_{2}$:
\flushpar
$(I)$ $H_{2} = \langle a \rangle \simeq C_{8}$:
\par
$(I_{1})$ $G \simeq C_{3} \times C_{8}$;
\par
$(I_{2})$ $G = \langle c, a \rangle \simeq C_{3} \rtimes C_{8}$, where 
$a^{-1}ca = c^{-1}$.
\flushpar
$(II)$ $H_{2} = \langle a , b \rangle \simeq C_{2} \oplus C_{4}$:
\par
$(II_{1})$ $G \simeq C_{3} \times (C_{2} \oplus C_{4})$;
\par
$(II_{2})$ $G = \langle c, a, b \rangle \simeq C_{3} \rtimes (C_{2} \oplus 
C_{4})$, where $a^{-1}ca = c$ and $b^{-1}cb = c^{-1}$.
\par
$(II_{3})$ $G = \langle c, a, b \rangle \simeq C_{3} \rtimes (C_{2} \oplus 
C_{4})$, where $a^{-1}ca = c^{-1}$ and $b^{-1}cb = c$.
\flushpar 
$(III)$ $H_{2} = \langle a_{1} , a_{2}, a_{3} \rangle \simeq C_{2}^{\oplus 3}$:
\par
$(III_{1})$ $G \simeq C_{3} \times C_{2}^{\oplus 3}$;
\par
$(III_{2})$ $G = \langle a_{1}, a_{2}, a_{3}, c \rangle \simeq 
C_{2}^{\oplus 3} \rtimes C_{3}$, where $c^{-1}a_{1}c = a_{1}$, $c^{-1}a_{2}c 
= a_{3}$ and $c^{-1}a_{3}c = a_{2}a_{3}$;
\par
$(III_{3})$ $G = \langle c, a, b \rangle \simeq C_{3} \rtimes 
C_{2}^{\oplus 3}$, where $a_{1}^{-1}ca_{1} = c$, $a_{2}^{-1}ca_{2} = c$ and 
$a_{2}^{-1}ca_{2} = c^{-1}$.
\flushpar 
$(IV)$ $H_{2} = \langle a, b \vert a^{4} = 1, a^{2} = b^{2}, 
b^{-1}ab = a^{-1} \rangle \simeq Q_{8}$:
\par
$(IV_{1})$ $G = \langle c \rangle \times 
\langle a, b \rangle \simeq C_{3} \times Q_{8}$;
\par
$(IV_{2})$ $G = \langle a, b, c \rangle \simeq Q_{8} \rtimes C_{3}$, where 
$c^{-1}ac = b$, $c^{-1}bc = ab$;
\par
$(IV_{3})$ $G = \langle c, a, b \rangle \simeq C_{3} \rtimes Q_{8}$, where 
$a^{-1}ca = c$, $b^{-1}cb = c^{-1}$.
\flushpar
$(V)$ $H_{2} = \langle a, b \vert a^{4} = 1, b^{2}=1, 
bab = a^{-1} \rangle \simeq D_{8}$:
\par
$(V_{1})$ $G = \langle c \rangle \times 
\langle a, b \rangle \simeq C_{3} \times D_{8}$;
\par
$(V_{2})$ $G = \langle c, a, b \rangle \simeq C_{3} \rtimes D_{8}$, where 
$a^{-1}ca = c$, $b^{-1}cb = c^{-1}$;
\par
$(V_{3})$ $G = \langle c, a, b \rangle \simeq C_{3} \rtimes D_{8}$, where 
$a^{-1}ca = c^{-1}$, $b^{-1}cb = c$;
\par
$(V_{4})$ $G \simeq S_{4}$. \qed
\endproclaim

As before, we denote by $A$ a target abelian threefold. In the case where 
$(I)$, $(II)$, $(III)$, $H_{2}$ is a commutative pre-C.Y. 
group of order $8$. However, this contradicts (2.5). In the case where 
$(IV_{1})$, $(IV_{3})$, $(V_{1})$, and $(V_{2})$, the subgroup 
$\langle a, c \rangle$ of $G$ is isomorphic to $C_{12}$. However, this 
again contradicts (2.5). In the case where $(V_{4})$, $G$ contains a subgroup 
which is isomorphic to $A_{4}$. However, this contradicts (2.14). Let us 
consider the case $(IV_{2})$. 
Set $H := \langle a, b \rangle$. Then, 
by (2.11), the representation $\rho_{H}$ of $H$ on $H^{0}(A, \Omega_{A}^{1})$ 
is decomposed as $\rho_{H} = \rho_{1, 0} \oplus \rho_{2,1}$. 
Let us write the $H$-stable subspace of $H^{0}(A, \Omega_{A}^{1})$ 
corresponding to $\rho_{1,0}$ by $V_{1}$. 
Then, $a(c(x)) = 
c(b(x)) = c(x)$ for $x \in V_{1}$ by $ac = cb$. Hence 
$V_{1}$ is $G$-stable. Therefore, by the Maschke Theorem, 
there exists a 2-dimensional $G$-stable subspace $V_{2}$ of 
$H^{0}(A, \Omega_{A}^{1})$ such that $H^{0}(A, \Omega_{A}^{1}) = V_{1} 
\oplus V_{2}$, and under an appropriate basis of $V_{1}$ and $V_{2}$, 
the matrix 
representation of $G$ on $H^{0}(A, \Omega_{A}^{1})$ is of the form; 
$a = \left(\matrix 1 & 0 & 0 \\ 0 & \zeta_{4} & 0 \\ 0 & 0 & -\zeta_{4} 
\endmatrix \right)$,  
$b = \left(\matrix 1 & 0 & 0 \\ 0 & 0 & \zeta_{4} \\ 0 & \zeta_{4} & 0 
\endmatrix \right)$ and $c = \left(\matrix \alpha & 0 \\ 0 & C \endmatrix 
\right)$, where $\alpha$ is a complex number and $C$ is a $2\times 2$ matrix. 
Since $c$ is of order $3$, $\alpha$ is either $1, \zeta_{3}$ or 
$\zeta_{3}^{-1}$. If $\alpha = 1$, then $H^{0}(A, \Omega_{A}^{1})^{G} = 
V_{1} \not= 0$. 
However, this contradicts our assumption that $G$ is a C.Y. group. Thus, 
we may assume that $\alpha = \zeta_{3}$ by replacing $c$ by $c^{-1}$ if 
necessary. Note that the eigen values of $C$ are now $\{1, \zeta_{3}^{-1}\}$.
Then the element $a^{2}c$ does not have an 
eigen value $1$, because $a^{2}c = \left(\matrix \zeta_{3} & 0 \\ 0 & -C 
\endmatrix \right)$, a contradiction to (2.4)(3). 
Hence the group in $(IV_{2})$ is 
not a C.Y. group. It remains to elminate the case $(V_{3})$. Set $V_{1} := 
H^{0}(A, \Omega_{A}^{1})^{c}$. Then, by (2.4)(3), $\text{dim}V_{1} = 1$. 
Using $ca = ac^{-1}$ and $cb = bc$, we see that $V_{1}$ is also 
$G$-stable. Then, again, by the 
Maschke Theorem, there exists a two-dimensional $G$-stable subspace $V_{2}$ 
of $H^{0}(A, \Omega_{A}^{1})$ such that $H^{0}(A, \Omega_{A}^{1}) = 
V_{1} \oplus V_{2}$. Note that by (2.11), this decomposition is also the 
irreducible decomposition of the representation of $\langle a, b \rangle 
(\simeq D_{8})$ and there exist basis of $V_{1}$ and $V_{2}$ 
under which we have $a = \left(\matrix 1 & 0 & 0 \\ 0 & \zeta_{4} & 0 
\\ 0 & 0 & -\zeta_{4} \endmatrix \right)$, 
$b = \left(\matrix -1 & 0 & 0 \\ 0 & 0 & 1 
\\ 0 & 1 & 0 \endmatrix \right)$ and $c = \left(\matrix 1 & 0 \\ 0 & C 
\endmatrix \right)$, where $C$ is a $2\times 2$ matrix. Then $bc$ is 
of the form $\left(\matrix -1 & 0 \\ 0 & D \endmatrix \right)$. Therefore, 
$\text{ord}(bc) = 2$ by (2.4)(3). On the other hand, 
since $bc = cb$, $\text{ord}(b) = 2$ and 
$\text{ord}(c) = 3$, we have $\text{ord}(bc) = 6$, 
a contradiction. Hence the group in $(V_{3})$ is not a C.Y. group, 
either.\qed 
\enddemo
This completes the proof (0.1) (I). \qed
\enddemo
\demo{Proof of (0.1) (II)}
\par
\vskip 4pt
Let us fix global coordinates $(z_{1}, z_{2}, z_{3})$ of $A$ as 
in (0.1)(I). 
Recall that under the identification 
$H^{*}(A, \Bbb C) = H_{\text{DR}}^{*}(A, \Bbb C)$
we have $H^{2}(A, \Bbb C) = \wedge^{2}H^{1}(A, \Bbb C)$ and  
$H^{1}(A, \Bbb Z) \otimes_{\Bbb Z}\Bbb C = H^{0}(A, \Omega_{A}^{1}) \oplus 
\overline{H^{0}(A, \Omega_{A}^{1})}$.  Using these identities and 
the description given in (0.1) (I) (1) and (2), we readily calculate that 
$$H^{2}(A, \Bbb C)^{G} 
= \Bbb C\langle dz_{1} \wedge d\overline{z}_{1}, dz_{2} \wedge d
\overline{z}_{2}, dz_{3} \wedge d\overline{z}_{3} \rangle \, \text{if}\, 
G \simeq C_{2}^{\oplus 2},$$ 
$$H^{2}(A, \Bbb C)^{G} = \Bbb C\langle dz_{1} 
\wedge d\overline{z}_{1}, dz_{2} \wedge d\overline{z}_{2} + dz_{3} \wedge 
d\overline{z}_{3} \rangle \, \text{if}\, G \simeq D_{8}.$$   
In addition, we have $c_{3}(X) = c_{3}(A)/\vert G \vert = 0$. 
Now the result follows from these equalities together with 
$c_{3}(X) = 2(\rho(X) - h^{1}(T_{X}))$ held for a Calabi-Yau 
threefold. \qed  
\enddemo
\demo{Proof of (0.1) (III)}
\par
\vskip 4pt 
We may construct explicit examples. 
\example{(2.17) Example for (I)(1) 
(Igusa's example [Ig, Page 678], [Ue, Example 16.16])} 
Let us first take three elliptic curves $E_{1}$, 
$E_{2}$ and $E_{3}$ and consider the product abelian 
threefold $A = E_{1} \times E_{2} \times E_{3}$. Let us fix points $\tau_{1} 
\in (E_{1})_{2} - \{0\}$, $\tau_{2} \in (E_{2})_{2} - \{0\}$,  
$\tau_{3} \in (E_{3})_{2} - \{0\}$ and define automorphisms $a$ and $b$ of 
$A$ by,
\par
$a(z_{1}, z_{2}, z_{3}) = (z_{1} + \tau_{1}, -z_{2}, -z_{3})$, and   
$b(z_{1}, z_{2}, z_{3}) = (-z_{1}, z_{2} + \tau_{2}, -z_{3} + \tau_{3})$.
\par 
Then, it is easy to check that 
$\langle a, b \rangle$ is isomorphic to $C_{2}^{\oplus 2}$ and 
acts on $A$ as a C. Y. group. Therefore, the quotient threefold 
$A/\langle a, b \rangle$ gives a desired example. \qed 
\endexample

\example{(2.18) Example for (I)(2)} Let us first take 
two elliptic curves $E_{1}$ and $E_{2}$ and consider the product abelian 
threefold $\tilde{A} = E_{1} \times E_{2} \times E_{2}$. 
Let us fix points $\tau_{1} \in (E_{1})_{4} - (E_{1})_{2}$, $\tau_{2}, 
\tau_{3} \in (E_{2})_{2}$ such that $\tau_{2} \not= \tau_{3}$ and define 
automorphisms $\tilde{a}$ and $\tilde{b}$ of $\tilde{A}$ by
\par 
$\tilde{a}(z_{1}, z_{2}, z_{3}) = (z_{1}+\tau_{1}, -z_{3}, z_{2})$, and   
$\tilde{b}(z_{1}, z_{2}, z_{3}) = (-z_{1}, z_{2}+\tau_{2}, 
-z_{3}+\tau_{3})$.
\par
Put $\tilde{G} = \langle \tilde{a}, \tilde{b} \rangle$. 
Then $\tilde{a}^{4} 
= \tilde{b}^{2} = id$, $\tilde{a}\tilde{b}\tilde{a}\tilde{b} = t_{\tau}$, 
$\tilde{a}t_{\tau}\tilde{a}^{-1} = t_{\tau}$ and 
$\tilde{b}t_{\tau}\tilde{b}^{-1} = t_{\tau}$, where $\tau = 
(0, \tau_{2}+\tau_{3}, \tau_{2}+\tau_{3})$. In 
particular, $\langle t_{\tau} \rangle (\simeq C_{2})$ is a normal subgroup of 
$\tilde{G}$. Set $A := \tilde{A}/\langle t_{\tau} \rangle$, 
$G := \tilde{G}/\langle t_{\tau} \rangle$, $a := \tilde{a}\, \text{mod}\, 
\langle t_{\tau} \rangle$, and $b := \tilde{b}\, \text{mod}\, 
\langle t_{\tau} \rangle$. 
Then $G = \langle a, b \rangle$ and acts on $A$ in the natural manner. 
It is now easy to check that this pair $(A, G)$, and hence, 
the quotient threefold $A/G$, gives a desired example. 
\qed 
\endexample
\enddemo 
 
Taking the contraposition of (II) and (III), we obtain the following 
criterion for the non-triviality of the second Chern class in terms of 
the Picard number:

\proclaim{Corollary (2.19)} Let $X$ be a Calabi-Yau threefold. 
If $\rho(X) = 1$ or $\rho(X) \geq 4$, then the second Chern class 
satisfies $c_{2}(X) \not= 0$ This is also optimal in the same sense 
as in (0.2).
\qed 
\endproclaim 
\demo{Proof of (0.1) (IV)}
\par
\vskip 4pt 
Let $A \rightarrow X$ be the minimal splitting covering of $X$ and 
$G = \langle a, b \rangle$ its Galois group as in (0.1)(I). 
Let us fix an origin 
$0 \in A$ and write $a = t_{\alpha} \circ a_{0}$ and 
$b = t_{\beta} \circ b_{0}$.  
We proceed our argument dividing into the cases (1) and (2) in (0.1) (I). 
\par
\vskip 4pt
{\it Case (1) in (0.1) (I).}  
\par
\vskip 4pt 
In this case $G \simeq C_{2}^{\oplus 2}$. 
Let us take the identity component $S_{1}$ of the kernel 
$0 \in S_{1} \subset \text{Ker}(a_{0} + id_{A} : A \rightarrow A)$, 
and consider the quotient map $\pi_{1} : A \rightarrow E_{1} := A/S_{1}$. 
This is an abelian fibration over an elliptic curve $E_{1}$ whose fibers are 
$S_{1} + p$, $p \in A$. Then, a similar argument to (2.13) shows that 
$p_{1} : A \rightarrow E_{1}$ is $G$-stable. Therefore, $\pi_{1} : A \rightarrow E_{1}$ induces an abelian fibration 
$\overline{\pi}_{1} : X = A/G \rightarrow E_{1}/G = \Bbb P^{1}$. 
Similarly, the identity components,  
$0 \in S_{2} \subset \text{Ker}(b_{0} + id_{A} : A \rightarrow A)$ and  
$0 \in S_{3} \subset \text{Ker}(a_{0}b_{0} + id_{A} : A \rightarrow A)$ 
induce abelian fibrations $\overline{\pi}_{2} : X  \rightarrow \Bbb P^{1}$ 
and $\overline{\pi}_{3} : X  \rightarrow \Bbb P^{1}$ respectively. Note that 
these three abelian fibrations $\overline{\pi}_{i}$ are mutually different 
by the shape of $a$ and $b$. Let us denote general fiber of 
$\overline{\pi}_{i}$ by $F_{i}$. 
\proclaim{Claim (2.20)} The classes $[F_{1}]$, $[F_{2}]$ and $[F_{3}]$ 
give a  basis of $\text{Pic}(X)_{\Bbb Q}$. 
\endproclaim 
\demo{Proof} Note that $F_{i} \cap F_{j} \not= \emptyset$ if 
$i \not= j$. Let $H$ be an ample divisor on $X$. 
Then, $(F_{i}.F_{j}.H) \not= 0$ for $i \not= j$ and  
$(F_{i}^{2}.H) = (F_{j}^{2}.H) = 0$. Therefore, $[F_{i}]$ and 
$[F_{j}]$ ($i \not= j$) are linearly independent in $\text{Pic}(X)_{\Bbb Q}$. 
Assume for a contradiction that $[F_{1}]$, $[F_{2}]$ and $[F_{3}]$ are 
linearly dependent in $\text{Pic}(X)_{\Bbb Q}$. 
Then there exist rational numbers $c_{1}$ and 
$c_{2}$ such that $F_{3} = c_{1}F_{1} + c_{2}F_{2}$ in $\text{Pic}(X)_{\Bbb Q}$ and satisfies 
$ 0 = (F_{3}^{2}.H) = 2c_{1}c_{2}(F_{1}.F_{2}.H)$. Therefore, either 
$c_{1} = 0$ or $c_{2} = 0$, a contradiction to the linear 
independence of $[F_{i}]$ and $[F_{j}]$ 
($i \not= j$). Since $\rho(X) = 3$, this gives the assertion. \qed 
\enddemo
\proclaim{Claim (2.21)} $\overline{\Cal A}(X) = \Bbb R_{\geq 0}[F_{1}] + 
\Bbb R_{\geq 0}[F_{2}] +  \Bbb R_{\geq 0}[F_{3}]$. 
\endproclaim 
\demo{Proof} The inclusion $\overline{\Cal A}(X) \supset 
\Bbb R_{\geq 0}[F_{1}] + \Bbb R_{\geq 0}[F_{2}] +  \Bbb R_{\geq 0}[F_{3}]$ 
is clear. Let us show the other inclusion. Let us choose an ample class 
$[H] \in \Cal A(X)$ 
and write $H = c_{1}F_{1} + c_{2}F_{2} + c_{3}F_{3}$ in $\text{Pic}(X)_{\Bbb R}$. Then $0 < (H.F_{1}.F_{2}) = c_{3}(F_{1}.F_{2}.F_{3})$. Since 
$(F_{1}.F_{2}.F_{3}) \geq 0$, we have $c_{3} > 0$. Similarly, $c_{1} > 0$ 
and $c_{2} > 0$. Therefore $\Cal A(X) \subset \Bbb R_{\geq 0}[F_{1}] + 
\Bbb R_{\geq 0}[F_{2}] +  \Bbb R_{\geq 0}[F_{3}]$. This implies the result. 
\qed 
\enddemo 
Since $h^{1}(\Cal O_{X}) = 0$, 
the set of numerically trivial classes of $\text{Pic}(X)$ is a 
finite group. Now combining this together with Claim (2.21) and the fact that 
$F_{i}$ are semi-ample, we also obtain the semi-ampleness assertion. \qed
\par
\vskip 4pt
{\it Case (2) in (0.1) (I).}
\par
\vskip 4pt
As before, take the identity components 
$0 \in S_{1} \subset \text{Ker}(a_{0}^{2} + id_{A}  : A \rightarrow A)$ and 
$0 \in E_{2} \subset \text{Ker}(a_{0} - id_{A}  : A \rightarrow A)$. 
Then, $S_{1}$ is an abelian surface and $E_{2}$ 
is an elliptic curve by the shape of $a_{0}$. Let us consider the quotient 
maps 
$\pi_{1} : A \rightarrow E := A/S_{1} $ and $\pi_{2} : A \rightarrow 
S := A/E_{2}$. Then as in (2.13), $G$ descends to the 
actions on the base spaces $E$ and $S$. Therefore $\pi_{1}$ and $\pi_{2}$ 
induce fibrations $\overline{\pi}_{1} : X \rightarrow \Bbb P^{1} = E/G$ and 
$\overline{\pi}_{2} : X \rightarrow \overline{S} = S/G$. Let $F_{1}$ be a 
general fiber of $\overline{\pi}_{1}$ and $F_{2}$ the pull back of an ample 
divisor on $\overline{S}$. Recall that in this case $\rho(X) = 2$ 
and $\partial \overline{\Cal A}(X)$ consists of two rays. Then, 
$[F_{1}]$ and $[F_{2}]$ are linearly independent in 
$\text{Pic}(X)_{\Bbb Q}$ and also satisfy
$\overline{\Cal A}(X) = \Bbb R_{\geq 0}[F_{1}] + \Bbb R_{\geq 0}[F_{2}]$. 
This implies the result. 
\qed 
\enddemo
\remark{Remark (2.22)} In Case (1), $X$ admits exactly 6 different non-trivial 
contractions corresponding to the 6 different strata of 
$\partial \overline{\Cal A}(X) - \{0\}$. 
Each of the three $1$-dimensional strata corresponds to an abelian fibration, 
as was observed in the proof, and each of the three $2$-dimensional strata 
corresponds to an elliptic fibration. Moreover, we can see that the base space 
of this elliptic fibration is a normal Enriques surface whose singularity is 
of Type $8A_{1}$ by using the shape of $a$ and $b$. 
In Case (2), we can also show that the base space $\overline{S}$ of 
$\overline{\pi}_{2} : X \rightarrow \overline{S}$ is a normal Enriques 
surface whose singularity is of Type $2A_{3} + 3A_{1}$. \qed
\endremark 

\head{\S 2K. Calabi-Yau threefolds of Type K}
\endhead
\par
\vskip 4pt
In this subsection, we study Calabi-Yau threefolds of Type K.  
The aim of this section is to show the following:  

\proclaim{Theorem (2.23)} Let $X$ be a Calabi-Yau threefold of Type K.  
Let $S \times E \rightarrow X$ be the minimal splitting cover, 
where $S$ is a K3 surface and $E$ is an elliptic curve, and $G$ its Galois 
group. Then, $X$ is isomorphic to $(S \times E)/G$ and, 
\flushpar
(I) $G$ is isomorphic to either $C_{2}^{\oplus n}$ 
($1 \leq n \leq 3$), $D_{2n}$ ($3 \leq n \leq 6$) or $C_{3}^{\oplus 2} 
\rtimes C_{2}$. 
\flushpar
(II) In each case, the Picard number $\rho(X)$, which is again equal to 
$h^{1}(T_{X})$, is 
uniquely determined by $G$ and is calculated as in the following table: 

$$\vbox{\tt \offinterlineskip
\def\vspa{height2pt & \omit && \omit && \omit && \omit && \omit && \omit 
    && \omit && \omit && \omit && \omit && \omit & \cr} 
\halign{&\vrule#&\strut~#\hfil~\cr 
\noalign{\hrule} 
\vspa
& \hfil $G$ && $C_{2}$ && $C_{2}^{\oplus 2}$ && $C_{2}^{\oplus 3}$ && $D_{6}$ 
&& $D_{8}$ && $D_{10}$ && $D_{12}$ && $C_{3}^{\oplus 2} \rtimes C_{2}$ & \cr
\vspa
\noalign{\hrule} 
\vspa
 & \hfil $\rho(X)$ && $11$  && $7$ && $5$ && $5$ && $4$ && $3$ && $3$ && $3$ 
& \cr
\vspa \noalign{\hrule} 
\hfil \cr }} $$ 
(III) The cases $G \simeq C_{2}^{\oplus n}$, 
where $1 \leq n \leq 3$, and $G \simeq D_{8}$ really occur. 
(Regrettably, it has not been settled yet 
whether the remaining cases occur or not.)
\flushpar
(IV) There exists a finite rational polyhedral cone $\Delta$ such that 
$\overline{\Cal A}'(X) = \text{Aut}(X) \cdot \Delta$, where 
$\overline{\Cal A}'(X)$ is the rational convex hull of the ample 
cone in $\text{Pic}(X)_{\Bbb R}$. Moreover, every nef $\Bbb Q$-divisor 
on $X$ is semi-ample. In particular, each $X$ admits only finitely many 
different contractions up to isomorphisms.  
\endproclaim 

\demo{Proof of (2.23)(I)} 
\enddemo
\par
\vskip 4pt
As in the subsection (2A), we define:

\definition{Definition (2.24)} We call a finite group $G$ a Calabi-Yau group 
of 
Type K, which again throughout this subsection, is abbriviated by a C.Y. 
group, if there exist a K3 surface $S$, an elliptic curve $E$ and a faithful 
representation $G \hookrightarrow \text{Aut}(S \times E)$ such that the 
following conditions (1) - (4) hold:
\roster
\item $G$ contains no elements of the form $(id_{S}, \text{non-zero 
translation of}\, E)$;
\item $g^{*}\omega_{S \times E} = \omega_{S \times E}$ for all $g \in G$;
\item $(S \times E)^{[G]} = \emptyset$;
\item $H^{0}(S \times E, \Omega_{S \times E}^{1})^{G} = \{0\}$. 
\endroster 
We call $S \times E$ a target threefold. \qed
\enddefinition 

As in subsection (2A), our proof is reduced to determine C. Y. group. 
\par
\vskip 4pt 

We first collect some easy Lemmas, whose varifications are so easy 
that we may omitt them. 

\proclaim{Lemma (2.25) ([Be2, Page 8, Proposition])} 
Let $S$ be a normal K3 surface and $E$ an elliptic curve. 
Then $\text{Aut}(S \times E) = \text{Aut}(S) \times \text{Aut}(E)$.  
In other word, each element $g$ of $\text{Aut}(S \times E)$ is 
of the form $(g_{S}, g_{E})$, where $g_{S} \in \text{Aut}(S)$ and 
$g_{E} \in \text{Aut}(E)$. \qed 
\endproclaim
\proclaim{Lemma (2.26)} Let $S$ be a K3 surface and $g$ an element of finite 
order of $\text{Aut}(S)$. Assume that $S^{[\langle g \rangle]} = \emptyset$. 
Then, 
\roster
\item 
if $g^{*}\omega_{S} = \omega_{S}$, then $g = id$; and 
\item 
if $g^{*}\omega_{S} \not= \omega_{S}$, then 
$g$ is of order $2$ and $g^{*}\omega_{S} = -\omega_{S}$. In this case, the 
quotient surface $S/\langle g \rangle$ is an Enriques surface. \qed
\endroster
\endproclaim
\proclaim{Lemma (2.27)}
Let $S$ and $E$ be same as in (2.25) and $G$ a finite 
subgroup of $\text{Aut}(S \times E)$ satisfying the same properties 
as (1), (2) and (3) in (2.24). Let 
$p_{1} : G \rightarrow \text{Aut}(S)$ and 
$p_{2} : G \rightarrow \text{Aut}(E)$ be 
the natural projections under the identification $\text{Aut}(S \times E) 
= \text{Aut}(S) \times \text{Aut}(E)$ (2.25) and put 
$G_{S} := \text{Im}(p_{1})$ 
and $G_{E} := \text{Im}(p_{2})$. Then, $G_{S} \simeq G \simeq G_{E}$ by 
$p_{1}$ and $p_{2}$. \qed 
\endproclaim
Let us return back to our study of C.Y. groups. 
\proclaim{Lemma (2.28)} Let $G$ be a C.Y. group and $S \times E$ its target 
threefold. Then, there exist a normal commutative subgroup $H$ of $G$ and 
an element $\iota$ of order $2$ of $G$ which satisfy the following properties 
(1) -(3):
\roster
\item 
$\iota \not\in H$ and $G = H \rtimes \langle \iota \rangle$, where the 
semi-direct product structure is given by $\iota h \iota = h^{-1}$ for all 
$h \in H$;
\item 
$\iota_{E} = -1_{E}$ and $H_{E} = \langle t_{a} \rangle \oplus \langle t_{b} 
\rangle$, under an appropriate origin of $E$, 
where $a$ and $b$ are torsion points such that 
$\text{ord}(a) \vert \text{ord}(b)$. In particular, $H \simeq H_{E} \simeq 
C_{n} \oplus C_{m}$ for some $1 \leq n \vert m$; and 
\item 
$S^{g_{S}} = \emptyset$ and $g_{S}^{*}\omega_{S} = - \omega_{S}$ 
for all $g \in G - H$, and $h_{S}^{*}\omega_{S} = \omega_{S}$ for all 
$h \in H$. 
\endroster
\endproclaim

\demo{Proof} 
Let us consider the natural homomorphism  
$G_{E} \rightarrow \text{GL}(H^{0}(E, \Omega_{E}^{1}))$ and denote by 
$H_{E}$ the kernel of this homomorphism. Then, we have  
$H^{0}(S \times E, \Omega_{S \times E}^{1})^{H} \simeq 
H^{0}(E, \Omega_{E}^{1})^{H_{E}} \simeq \Bbb C$.
In particular, $H_{E} \not= G_{E}$. Take an arbitrary element $\iota_{E}$  
of $G_{E} - H_{E}$ and put 
$\iota := (\iota_{S}, \iota_{E}) \in \text{Aut}(S \times E)$. 
Then, there exists a complex number $\alpha \not= 1$ such that 
$\iota_{E}^{*}\omega_{E} = \alpha\omega_{E}$. 
This implies $E^{\iota_{E}} \not= \emptyset$ and  
$\iota_{S}\omega_{S} = \alpha^{-1}\omega_{S}$. In particular, 
$S^{\iota_{S}} = \emptyset$. Therefore $\iota_{S}$ is an involution 
and $\alpha = -1$ by (2.26). Let us fix one of such an $\iota$. 
Then, for any $\iota_{E}' \in G_{E} - H_{E}$, we have 
$\iota_{E}' \circ \iota_{E} \in H_{E}$. Therefore,
$G_{E} = H_{E} \rtimes \langle \iota_{E} \rangle$. 
Fix the origin $0$ in $E^{\iota_{E}}$. 
Then $\iota_{E} = -1_{E}$ and $-1_{E} \circ t_{a} \circ -1_{E} = t_{-a} 
= t_{a}^{-1}$. In particular, $\iota_{E} \circ h_{E} \circ \iota_{E} 
= h_{E}^{-1}$ if $h \in H$. 
This gives the semi-direct product structure. 
Moreover, since $H_{E}$ consists of translations of $E$, there exist
positive integers $n$ and $m$ such that $H_{E} \simeq C_{n} \oplus C_{m}$ 
and that $n \vert m$. 
\qed 
\enddemo 
\proclaim{Lemma (2.29)} Let $(n, m)$ be same as in (2.28)(2). Then, $(n, m) 
\in \{(1, k) (1 \leq k \leq 6), (2, 2), (3, 3)\}$. 
\endproclaim

\demo{Proof} For proof, we make use of the following result due to Nikulin:

\proclaim{Theorem (2.30) ([Ni, Page 106, Section 5, Paragraph 8])} 
Let $S$ be a K3 surface. 
\roster
\item 
Let $g \not= id$ be a Gorenstein automorphism of finite order. 
Then, $\text{ord}(g) \leq 8$. 
Moreover, $S^{g}$ is a finite set and its cardinality $\vert S^{g} \vert$ 
is given as in the following table:
$$\vbox{\tt \offinterlineskip
\def\vspa{height2pt & \omit && \omit && \omit && \omit && \omit && \omit 
    && \omit && \omit && \omit && \omit && \omit & \cr} 
\halign{&\vrule#&\strut~#\hfil~\cr 
\noalign{\hrule} 
\vspa
& \hfil $\text{ord}(g)$ && $2$ && $3$ && $4$ && $5$ && $6$ && $7$ && $8$ & \cr
\vspa
\noalign{\hrule} 
\vspa
 & \hfil $\vert S^{g} \vert$ && $8$  && $6$ && $4$ && $4$ && $2$ && $3$ && 
$2$ & \cr
\vspa \noalign{\hrule} 
\hfil \cr }}.$$ 
\item 
Let $H$ be a finite, commutative, Gorenstein subgroup of $\text{Aut}(S)$. 
Then $H$ is isomorphic to either one of 
$C_{k}$ ($1 \leq k \leq 8$), $C_{2}^{\oplus l}$ ($2 \leq l \leq 4$), 
$C_{2} \oplus C_{4}$, $C_{2} \oplus C_{6}$, $C_{3}^{\oplus 2}$, or 
$C_{4}^{\oplus 2}$. \qed 
\endroster
\endproclaim
 
Recall that $G_{S} \simeq G_{E}$ and $H_{S} \simeq H_{E}$ under 
the isomorphisms in (2.27). Then, $H_{S}$ is a Gorenstein 
automorphism group of $S$ and is isomorphic to $H_{S} \simeq C_{n} \oplus 
C_{m}$. Now it is sufficient to elminate the following cases in (2.30)(2): 
$$(n, m) = (1,7),\, (1,8),\, (2,4),\, (2,6),\, (4,4).$$ 
\flushpar
\vskip 4pt 
We elminate all cases by more or less similar method. So, we explain 
how to do this for the hardest case $(n, m) = (2,4)$ and leave the other 
cases to the readers.  Assume to the contrary, that $(n, m) = (2, 4)$. 
Then $H_{S} = \langle g_{S} 
\rangle \oplus \langle h_{S}  \rangle \simeq C_{2} \oplus C_{4}$.
Note that $\langle g_{S}, h_{S}, \iota_{S} \rangle / \langle h_{S}^{2} 
\rangle$ is isomorphic to $C_{2}^{\oplus 3}$ and acts on 
$S^{h_{S}^{2}} - S^{h_{S}}$. Note also that by (2.33)(1), 
we have $\vert S^{h_{S}^{2}} - S^{h_{S}} \vert = 4$. 
Then, this action induces a homomorphism $\varphi : C_{2}^{\oplus 3} 
\rightarrow S_{4}$. Since $S_{4}$ does not contain a subgroup isomorphic to 
$C_{2}^{\oplus 3}$, 
we have $\text{Ker}(\varphi) \not= \{id\}$. Moreover, 
$\text{Ker}(\varphi) \subset 
\langle g_{S}, h_{S} \rangle / \langle h_{S}^{2} \rangle$,
because 
$S^{f} = \emptyset$ for $f \in G_{S} -H_{S}$.
Let $\alpha \in \langle g_{S}, h_{S} \rangle$ be a lift of a non-trivial 
element of $\text{Ker}(\varphi)$ and take $P \in S^{h_{S}^{2}} - S^{h_{S}}$. 
Then we have a natural 
injection $\langle \alpha, h_{S}^{2} \rangle \hookrightarrow 
\text{SL}(T_{S, P}) = \text{SL}(2, \Bbb C)$. 
In addition, using $h_{S} \not\in \text{Ker}(\varphi)$, 
we see that $\langle \alpha, h_{S}^{2} \rangle$ is isomorphic to either 
$C_{2}^{\oplus 2}$ or $C_{2} \oplus C_{4}$. However, this 
contradicts the following well-known:
\proclaim{Theorem (2.31) (see for example [Su, Chap. III, \S 6])}
Let $G$ be a finite subgroup of $\text{SL}(2, \Bbb C)$. 
Then $G$ is isomorphic to either one of 
$C_{n}$, $Q_{4n}$, $T_{24}$, $O_{48}$ or $I_{120}$, where $T_{24}$, 
$O_{48}$, $I_{120}$ are the binary polyhedral groups of indicated orders. 
\qed 
\endproclaim
This completes the proof of (2.23)(I). \qed
\enddemo
\demo{Proof of (2.23)(II)}
\par
\vskip 4pt
The equality $\rho(X) = h^{1}(T_{X})$ follows from the same 
reason as for Type A. Since $\rho(X) = \text{dim} H^{2}(S \times E, 
\Bbb C)^{G} = \text{dim}H^{2}(S, \Bbb C)^{G_{S}} + 1$, 
the proof is reduced to the calculation of 
$\text{dim}H^{2}(S, \Bbb C)^{G_{S}}$ for each $G$. Our calculation is based 
on the topological Lefschetz fixed point formula, (2.28)(3) and (2.30)(1) and 
is similar for all $G$. So, we explain how 
to calculate $\rho(X)$ only for $G \simeq D_{8}$ and leave the other cases 
to the reader. Write $G_{S} = \langle 
a , b \vert a^{4} = b^{2} = 1, bab = a^{-1} \rangle (\simeq D_{8})$ and 
write the irreducible decomposition of $G \curvearrowright H^{2}(S, \Bbb C)$ 
as $H^{2}(S, \Bbb C) = \rho_{1, 0}^{\oplus p} \oplus \rho_{1, 1}^{\oplus q} 
\oplus \rho_{1, 2}^{\oplus r} \oplus \rho_{1, 3}^{\oplus s} \oplus 
\rho_{2, 1}^{\oplus t}$, where we adopt the same notation as in (2.10). 
Then, by $\text{dim} H^{2}(S, \Bbb C) = 22$, we have 
$22 = p + q + r + s + 2t$. 
Using $\vert S^{a} \vert = 4$ (2.30)(1), and applying the topological 
Lefschetz formula, we get
$4 = \vert S^{a} \vert = 2 + \text{tr}(a^{*} \vert H^{2}(S, \Bbb C)) 
= 2 + p + q - r - s$, that is, $2 = p + q - r - s$.  
Similarly, from $\vert S^{a^{2}} \vert = 8$ by (2.30)(1), $\vert S^{b} 
\vert = 0$ and $\vert S^{ab} \vert = 0$ by (2.28)(3), we obtain
$6 = p + q + r + s - 2t$, $-2 = p - q + r - s$, $-2 = p - q - r + s$. 
Now solving this system of equations, we raedily find that 
$p = 3$, $q = 5$, $r = s = 3$ and $t = 4$. 
Therefore, $\text{dim}H^{2}(S, \Bbb C)^{G_{S}} = p = 3$.  
Hence, $\rho(X) = 3 + 1 = 4$. \qed 
\enddemo 

\demo{Proof of (2.23)(III)}
\par
\vskip 4pt  
We may construct a Calabi-Yau threefold of Type K such that the Galois 
group $G$ of its minimal splitting cover is 
isomorphic to $C_{2}^{\oplus n}$ $(1 \leq n \leq 3)$ 
and $D_{8}$ respectively.
\enddemo
\example{(2.32) Example for $C_{2}^{\oplus n}$ $(1 \leq n \leq 3)$} 
Let us first take three elliptic curves (with fixed origins) 
$E_{1}$, $E_{2}$ and $E$ and denote by $S := \text{Km}(E_{1} \times E_{2})$ 
the smooth Kummer surface associated with the product abelian surface $E_{1} 
\times E_{2}$. Choose elements $a_{i}, b_{i} \in (E_{i})_{2} - \{0\}$ 
such that $a_{i} \not= b_{i}$ for each $i = 1, 2$. Then, the following three 
automorphisms of $E_{1} \times E_{2}$ descend to those of $\text{Aut}(S)$: 
\par
$(z_{1}, z_{2}) \mapsto (-z_{1} + a_{1}, -z_{2} +a_{2})$, \, 
$(z_{1}, z_{2}) \mapsto (z_{1} + b_{1}, z_{2})$, \,  
$(z_{1}, z_{2}) \mapsto (z_{1}, z_{2} + b_{2})$. 
\par 
We denote them by $\theta$, $t_{1}$ and $t_{2}$ respectively. Let us 
choose $P_{1}, P_{2} \in 
(E)_{2} -\{0\}$ such that $P_{1} \not= P_{2}$ 
and define the three automorphisms of 
$S \times E$ by $\tilde{\theta} := (\theta, -1_{E})$, 
$\tilde{t_{1}} := (t_{1}, t_{P_{1}})$ and 
$\tilde{t_{2}} := (t_{2}, t_{P_{2}})$. 
Set $G_{1} := \langle \tilde{\theta} \rangle$, 
$G_{2} := \langle \tilde{\theta}, \tilde{t_{1}} \rangle$ and 
$G_{3} := \langle \tilde{\theta}, \tilde{t_{1}}, \tilde{t_{2}} \rangle$. 
Then, $G_{n} \simeq C_{2}^{\oplus n}$ and act on $S \times E$ 
as C. Y. groups. Therefore, the quotient threefolds 
$(S \times E)/G_{n}$ give desired examples. \qed 
\endexample 

Next we construct an explicit example for $D_{8}$. Let us first observe 
the following: 

\proclaim{Proposition (2.33)} 
Let $D_{8} = \langle a, b \vert a^{4} = b^{2} = 1, bab = a^{-1} \rangle$ 
and $V$ the regular representation of $D_{8}$ defined by 
$V = \rho_{2,1} \oplus \rho_{1, 0} \oplus \rho_{1,2} \oplus \rho_{1, 1} 
\oplus \rho_{1,3}$.
Regard $\Bbb P^{5} = \text{Proj}(\oplus \text{Sym}^{\cdot}V)$ and 
define $S$ to be the complete intersection in $\Bbb P^{5}$ given by 
\par
$x_{0}^{2} + x_{1}^{2} + x_{2}x_{3} + x_{4}x_{5} = x_{0}x_{1} + x_{2}^{2} 
+ x_{3}^{2} + x_{4}^{2}  + x_{5}^{2} = Ax_{2}x_{3} + Bx_{4}x_{5} = 0,$
\par
where $A$ and $B$ are sufficiently general complex numbers. Then, 
$S$ is a K3 surface and is stable under the action of $D_{8}$ on $\Bbb P^{5}$. 
Moreover, the induced action on $S$ is faithful 
and satisfies $a^{*}\omega_{S} = \omega_{S}$, 
$b^{*}\omega_{S} = -\omega_{S}$ and 
$S^{a^{k}b} = \emptyset$ for $k = 0, 1, 2, 3$. 
\endproclaim 
\demo{Proof} Smoothness of $S$ follows from the Jacobian criterion. 
The rest follows from direct calculations. 
 \qed \enddemo 

\example{(2.34) Example for $D_{8}$} 
Let $S$ be a K3 surface in (2.33) and $E$ be an elliptic curve. 
Let us consider elements of $\text{Aut}(S \times E)$ defined by 
$g := (a, t)$ and $\iota := (b, -1_{E})$, where $a$ and $b$ are same as 
in (2.33) and $t$ is a translation automorphism of $E$ of order 4. 
Then, $\langle g, \iota \rangle$ is isomorphic to $D_{8}$ and acts on 
$S \times E$ as a C. Y. group. Therefore, $(S \times E)/\langle g, \iota 
\rangle$ gives a desired example. \qed
\endexample 

\demo{Proof of (2.23)(IV)}
\par
\vskip 4pt  
We may identify 
$\text{Pic}(X)_{\Bbb Q} = (\text{Pic}(S \times E)^{G})_{\Bbb Q}$ 
via the quotient map. By using $h^{1}(\Cal O_{S}) = 0$, (2.25) and 
the Kunneth formula, 
we also obtain $(\text{Pic}(S \times E)^{G})_{\Bbb Q} 
= (\text{Pic}(S)^{G_{S}})_{\Bbb Q} \oplus \Bbb Q$. Now the result follows 
from (1.9) and the semi-ampleness of rational nef divisor 
on K3 surface, because again the torsion group of $\text{Pic}(X)$ is finite 
by $h^{1}(\Cal O_{X}) = 0$. 
\qed
\enddemo 

As an immediate Corollary of (0.1) and (2.23), we obtain Corollary (0.2) 
in Introduction and the following:

\proclaim{Corollary (2.35)} Let $X$ be a Calabi-Yau threefold. Assume that 
$\pi_{1}(X)$ is infinite. Then $\pi_{1}(X)$ falls into one of the following 
exact sequences:
\par
$0 \rightarrow  \Bbb Z^{\oplus 6} \rightarrow \pi_{1}(X) \rightarrow G 
\rightarrow 1$, where $G$ is isomorphic to either $C_{2}^{\oplus 2}$ or 
$D_{8}$ or,
\par
$0 \rightarrow  \Bbb Z^{\oplus 2} \rightarrow \pi_{1}(X) \rightarrow G 
\rightarrow 1$, where $G$ is isomorphic to either $C_{2}^{\oplus n}$ 
($1 \leq n \leq 3$), $D_{2n}$ ($3 \leq n \leq 6$) or $C_{3}^{\oplus 2} 
\rtimes C_{2}$. In particular, $\pi_{1}(X)$ is always solvable. \qed
\endproclaim
\head
{\S 3. Classification of $c_{2}$-contractions of Calabi-Yau threefolds}
\endhead
\par
\vskip 4pt 
In this section, we give a generalisation (and a correction) 
of our earlier work of $c_{2} = 0$ contractions of simply connected 
Calabi-Yau threefolds [Og 1-4]. 
First we remark the following easy facts on abelian varieties 
applied in Sections 3 and 4.  
 
\proclaim{Lemma (3.1)} Every contraction of an abelian variety $A$ 
is of the form of the exact sequence of abelian varieties: 
$0 \rightarrow F \rightarrow A \rightarrow \overline{A} \rightarrow 0$.
\endproclaim 
\demo{Proof} Let $f : A \rightarrow \overline{A}$ be a contraction 
and $F$ a smooth fiber. 
Let us choose an origin $0$ of $A$ in $F$ and 
regard $A$ as a group variety. Since $[t_{-a}(F)] = [F] \in 
H^{*}(A, \Bbb Z)$ for all $a \in A$, we see that $f(t_{-a}(F))$ 
is also a point by taking an appropriate intersection with the 
pull back of an ample divisor on $\overline{A}$. Therefore, 
$0 \in t_{-a}(F) \cap F$ and $t_{-a}(F) = F$ for all $a \in F$.  
Hence, $F$ is an abelian subvariety of $A$ and induces an isomorphism
$\overline{A} \simeq A/F$. This implies the result. \qed
\enddemo 
\proclaim{Proposition (3.2)} 
\roster
\item 
Let $(A, h)$ be a pair of an abelian variety of $\text{dim}A = n$ and 
its automorphism $h$ such that $h^{*} \vert H^{0}(A, \Omega_{A}^{1}) 
= \zeta_{3}$, the scalar multiplication by $\zeta_{3}$. Then 
$(A, g)$ is isomorphic to the pair $(E_{\zeta_{3}}^{n}, 
\text{diag}(\zeta_{3}, \cdots , \zeta_{3}))$. 
\item 
Any $\text{diag}(\zeta_{3}, \cdots , \zeta_{3})$-satble contraction of 
$E_{\zeta_{3}}^{n}$ is $\text{diag}(\zeta_{3}, \cdots , \zeta_{3})$-
equivariantly isomorphic to the projection 
$p_{1, \cdots , m} : E_{\zeta_{3}}^{n} \rightarrow E_{\zeta_{3}}^{m}$ 
to the first $m$-factors for some $m$. 
\endroster
\endproclaim 
\demo{Proof} The assertion (1) is shown in [CC, Proposition 5.7] and also 
follows from the argument of [Og3, Section 1]. Let us show the assertion (2). 
For the sake of simplicity, we put $A = E_{\zeta_{3}}^{n}$ and 
$g = \text{diag}(\zeta_{3}, \cdots , \zeta_{3})$. We regard 
the universal cover $\Bbb C^{n}$ as a $\Bbb Z[\zeta_{3}]$-module via the 
scalar action of $g$. Write  $A = \Bbb C^{n}/\Lambda$. Then, the lattice 
$\Lambda$ is a $\Bbb Z[\zeta_{3}]$-submodule 
of $\Bbb C^{n}$ and coincides with the subset 
$\Bbb Z[\zeta_{3}]^{\oplus n} \subset \Bbb C^{n}$. 
Let $\varphi : A \rightarrow B$ a $g$-stable contraction of $A$ and 
take the fiber $F$ of $\varphi$ which contains the origin $0 \in A$. 
Then, $F$ is an abelian subvariety of $A$ by (3.1). 
Let us denote by $\Lambda_{F}$ the sublattice of $\Lambda_{A}$ corresponding 
to $F$. Since $F$ is $g$-stable, $\Lambda_{F}$ is also a 
$\Bbb Z[\zeta_{3}]$-submodule of $\Lambda_{A}$ of rank $n - m$, where 
$m = \text{dim}(B)$. Moreover, $\Lambda_{F}$ is primitive in $\Lambda_{A}$, 
because $\Lambda_{F} = \Lambda_{A} \cap V_{F}$, where $V_{F}$ is the linear 
subspace of $\Bbb C^{n}$ corresponding to $F$. Therefore, there exists an 
element 
$h \in \text{GL}(n, \Bbb Z[\zeta_{3}])$, where we regard 
$\text{GL}(n, \Bbb Z[\zeta_{3}])$ as a subgroup of 
$\text{Aut}(\Lambda_{A})$,  
such that $h(\Lambda_{F}) = \{0\} \oplus \Bbb Z[\zeta_{3}]^{\oplus (n-m)}$. 
Recall that $\Bbb Z[\zeta_{3}]$ is an Euclidean domain and is 
also the endmorphism ring of $E_{\zeta_{3}}$. Then, by the elementary divisor 
theory, and by the fact that $g$ is contained in the center of 
$\text{Aut}_{\text{Lie}}(A)$, we see that the image 
of the natural representation $\text{Aut}_{\text{Lie}}(E_{\zeta_{3}}^{n}) 
\rightarrow \text{Aut}(\Lambda_{A})$ coincides with 
$\text{GL}(n, \Bbb Z[\zeta_{3}])$. In particular, $h$ is the image of some 
Lie automorphism $\tilde{h}$ of $A$. It is clear that this $\tilde{h}$ gives 
a desired $g$-equivariant isomorphism. 
\qed 
\enddemo

\proclaim{Theorem (3.3)} 
Let $\Phi : X \rightarrow W$ be a $c_{2}$-contraction. Assume that $\Phi$ 
is an isomorphism. Then $X$ is a smooth Calabi-Yau threefold of Type A (0.1). 
\endproclaim 
\demo{Proof} Since $(c_{2}(X).H) = 0$ for ample divisors on $X$, 
we see that $c_{2}(X) = 0$ as a linear form on 
$\text{Pic}(X)_{\Bbb R}$. 
\qed
\enddemo
Next we consider non-trivial birational $c_{2}$-contraction. 

\proclaim{Theorem (3.4) (cf. [Og3, Main Theorem])} 
Let $\Phi : X \rightarrow W$ be a $c_{2}$-contraction.
Assume that $\Phi$ is birational but not an isomorphism. Then, 
$\Phi : X \rightarrow W$ is isomorphic to either one of the following:
\roster
\item"(1)" 
The unique crepant resolution $\Phi_{7} : X_{7} \rightarrow 
\overline{X}_{7} := A_{7}/\langle g_{7} \rangle$ of $\overline{X}_{7}$, 
where $(A_{7}, g_{7})$ is the Klein pair. In this case $\rho(X_{7}) = 24$, 
$\rho(\overline{X}_{7}) = 3$ and $\pi_{1}(X_{7}) = \{1\}$. 
\item"(2-0)" 
The unique crepant resolution $\Phi_{3} : X_{3} \rightarrow 
\overline{X}_{3} := A_{3}/\langle g_{3} \rangle$, 
where $(A_{3}, g_{3})$ is the Calabi pair. 
In this case $\rho(X_{3}) = 36$, 
$\rho(\overline{X}_{3}) = 9$ and $\pi_{1}(X_{3}) = \{1\}$. 
\item"(2-1)" 
The unique crepant resolution $\Phi_{3,1} : X_{3,1} \rightarrow 
\overline{X}_{3,1} := A_{3}/\langle g_{3}, h \rangle$
of $\overline{X}_{3,1}$, 
where $(A_{3}, g_{3})$ is the Calabi pair and 
$\langle g_{3}, h \rangle \simeq C_{3}^{\oplus 2}$. 
Moreover, $\langle h \rangle$ acts on $\overline{X}_{3}$ freely and the 
representation of $\langle g_{3}, h \rangle$ on $H^{0}(\Omega_{A_{3}}^{1})$ 
is given by:
\par 
$g_{3} \mapsto 
\left(\matrix \zeta_{3} & 0 & 0 \\ 0 & \zeta_{3} & 0 \\ 
0 & 0 & \zeta_{3}\endmatrix \right)$ 
and $h \mapsto \left(\matrix 1 & 0 & 0 \\ 0 & \zeta_{3} & 0 \\ 
0 & 0 & \zeta_{3}^{2}\endmatrix \right)$. 
\par 
In this case $\rho(X_{3,1}) = 12$, $\rho(\overline{X}_{3,1}) = 3$ and 
$\pi_{1}(X_{3}) \simeq C_{3}$.
\item"(2-2)"
The unique crepant resolution $\Phi_{3,2} : X_{3,2} \rightarrow 
\overline{X}_{3,2} := A_{3}/\langle g_{3}, h, k \rangle$, 
where $(A_{3}, g_{3})$ is again the Calabi pair and  
$\langle g_{3}, h, k \rangle$ is 
the unique non-commutative group of order 27 whose elements 
($\not=$) are all of order 3 (cf. [Bu, Chap.8, Page 158]). 
Moreover, $\langle g_{3}, h, k\rangle /\langle g_{3} \rangle$ 
is isomorphic to $C_{3}^{\oplus 2}$ and acts on $\overline{X}_{3}$ freely 
and the representaion of $\langle g_{3}, h, k \rangle$ on 
$H^{0}(\Omega_{A_{3}}^{1})$ is given by:
\par 
$g_{3} \mapsto \left(\matrix \zeta_{3} & 0 & 0 \\ 0 & \zeta_{3} & 0 \\ 
0 & 0 & \zeta_{3}\endmatrix \right)$,  
$h \mapsto \left(\matrix 1 & 0 & 0 \\ 0 & \zeta_{3} & 0 \\ 
0 & 0 & \zeta_{3}^{2}\endmatrix \right)$, and 
$k \mapsto \left(\matrix 0 & 0 & 1 \\ 1 & 0 & 0 \\ 
0 & 1 & 0\endmatrix \right).$
\par
In this case $\rho(X_{3,2}) = 4$, 
$\rho(\overline{X}_{3,2}) = 1$ and 
$\pi_{1}(X_{3,2}) \simeq C_{3}^{\oplus 2}$.
\endroster
\endproclaim
\demo{Proof} This is a refinement of [Og3], in which we have 
already obtained the following properties based on [SBW, Main Theorem] and 
the PID property of the cyclotomic integer rings $\Bbb Z[\zeta_{I}]$ 
of relatively small degree $\varphi(I)$ ([MM, Main Theorem]): 
\proclaim{Lemma (3.5) ([Og3, Key Claim Page 334, Lemma (2.1), (2.2), 
Remark after Theorem 3])} 
Under the same assumption of (3.4), $\Phi : X \rightarrow W$ is 
isomorphic to either $\Phi_{7} : X_{7} \rightarrow 
\overline{X}_{7}$ or the (necessarily unique) crepant 
resolution of the quotient $A_{3}/G$, where $G$ satisfies that:
\roster
\item"(i)"
$G$ is a finite Gorenstein automorphism group of $A_{3}$ and 
$\langle g_{3} \rangle \vartriangleleft G$;
\item"(ii)"
$G$ contains no non-trivial translations; and that,
\item"(iii)"
the induced action of $G/\langle g_{3} \rangle$ on 
$A_{3}/\langle g_{3} \rangle$ 
is fixed point free. \qed
\endroster
\endproclaim 
Therefore, it is sufficient to determine such $G$ of. 
Write $A := A_{3}$. Then, by (i) and (ii), 
the natural homomorphism $G \rightarrow \text{SL}(H^{0}(A, \Omega_{A}^{1}))$ 
is injective. In particular, $g_{3}$ is contained in the center of $G$. 
Take $h \in G - \langle g_{3} \rangle$ and set $d := \text{ord}(h)$. 
Then, $d$ is either $3$, $9$ or $27$, because $h$ acts on the set 
$A^{g_{3}}$ freely and $\vert A^{g_{3}} \vert = 27$.  Moreover, at least one 
eigen value of $h^{*} \vert H^{0}(A, \Omega_{A}^{1})$ must be one, 
because, otherwise $h$ has an isolated fixed points. But this contradicts 
$h \in G - \langle g_{3} \rangle$ if $d = 3$ and the condition $(iii)$ if 
$d \not= 3$. Now, repeating the same argument as in (2.4)(3), 
we see that $\varphi(d) \leq 2 = (6-2)/2$. Therefore, $d = 3$ 
and $\vert G \vert = 3^{n}$ for some positive integer $n$. 
Let us consider a maximal normal commutative subgroup $H$ of $G$ and put 
$\vert H \vert = 3^{m}$. Then $g_{3} \in H$ and $m(m+1) \geq 2n$ by (2.7). 
Assume that $n \geq 4$. Then $m \geq 3$. Therefore, 
$H$ contains a subgroup $L$
such that $L \simeq C_{3}^{\oplus 3}$ and $g_{3} \in L$. 
Set $L = \langle g_{3} \rangle \oplus 
\langle h \rangle \oplus \langle k \rangle$. Then, there exists a basis 
of $H^{0}(A, \Omega_{A}^{1})$ under which (after replacing $k$ by $k^{2}$ if 
necessarily) the matrix representation of $L$ is of the form:
$g_{3}^{*} = \text{diag}(\zeta_{3},\zeta_{3}, \zeta_{3})$, 
$h^{*} = \text{diag}(1,\zeta_{3}, \zeta_{3}^{2})$, 
and 
$k^{*} = \text{diag}(\zeta_{3}, 1, \zeta_{3}^{2})$ or 
$\text{diag}(\zeta_{3}, \zeta_{3}^{2}, 1)$. 
However, this implies either $g_{3}^{2}h = k$ or $g_{3}h = k$, 
a contradiction. 
Therefore $n \leq 3$ and $G \not\simeq C_{3}^{\oplus 3}$. This together 
with conditions (i) - (iii) readily implies the result. \qed
\enddemo
Let us next consider $c_{2}$-contractions $\Phi : X \rightarrow W$ 
such that $\text{dim}(W) = 2$. It is known by [Nk, Corollary (0.4)] 
(see also [OP, Corollary (2.6)]) that $(W, 0)$ is klt and that 
$\Cal O_{W}(12K_{W}) \simeq \Cal O_{W}$. 
Let us define the global canonical index of $W$ by 
$I := I(W) := \text{min}\{n \in \Bbb Z_{> 0} \vert 
\Cal O_{W}(nK_{W}) \simeq \Cal O_{W}\}$ and take 
the global index one cover of $W$: 
$$\pi : T := Spec_{\Cal O_{W}}(\oplus_{k=0}^{I-1} 
\Cal O_{W}(-kK_{W})) \rightarrow W.$$ 
Then, by [Kaw1, Pages 608-609] (see also [Zh, \S 2]),  
$T$ is either a normal K3 surface or 
a smooth abelian surface and that $\pi$ is a cyclic Galois covering 
\'etale in codimension one such that the Galois group $G$, which is isomorphic 
to $C_{I}$, acts faithfully on 
$H^{0}(\Cal O_{T}(K_{T})) = \Bbb C \omega_{T}$. Note that 
$I \vert 12$ and $I \geq 2$. (Indeed, if $I = 1$, then $W$ is either a 
normal K3 surface or an abelian surface. However, 
this contradicts $h^{2}(\Cal O_{X}) = h^{1}(\Cal O_{X}) = 0$.)
We call $\Phi : X \rightarrow W$ of Type IIA if $T$ 
is a smooth abelian surface and of Type IIK if $T$ is a normal K3 surface. 
It has been shown in [Og2] the following: 
\proclaim{Theorem (3.6) ([Og2, Main Theorem])} Let $\Phi : X \rightarrow W$ 
be a $c_{2}$-contraction of Type IIA. 
Then $\Phi : X \rightarrow W$ is isomorphic to one of the 
relatively minimal models over $E_{\zeta_{3}}^{2}/\langle 
\overline{g}_{3} \rangle$ of 
$$p_{12} : X_{3} \overset\Phi_{3}\to\longrightarrow 
\overline{X}_{3} = E_{\zeta_{3}}^{3}/\langle g_{3} 
\rangle 
\overset\overline{p}_{12}\to\longrightarrow E_{\zeta_{3}}^{2}/\langle 
\overline{g}_{3} \rangle,$$ 
where $\Phi_{3} : X_{3} \rightarrow \overline{X}_{3}$ 
is the contraction in (3.4)(2-0) and $\overline{p}_{12}$ is the natural 
projection. In particular, $X$ is smooth 
and $I(W) = 3$. Moreover, there exist exactly $2^{9}$ 
such relatively minimal models. 
\qed
\endproclaim   
The next Theorem is a generalisation and also a correction of [Og4]: 
\proclaim{Theorem (3.7) (cf. [Og4, Main Theorem])} Let 
$\Phi : X \rightarrow W$ be a $c_{2}$-contraction of Type IIK. Then 
$\Phi : X \rightarrow W$ is isomorphic to 
either: 
\roster
\item 
A fiber space of a Calabi-Yau threefold of Type A in (0.1)(I)(1) 
corresponding to a $2$-dimensional face of its nef cone. In this case, 
$I(W) = 2$ and $\rho(W) = 2$. (See also (2.22).)
\item 
The fiber space of a Calabi-Yau threefold of Type A in (0.1)(I)(2) 
given by the boundary of its nef cone corresponding to the 
elliptic fibration. In this case, $I(W) = 2$ and $\rho(W) = 1$. 
(See also (2.22).)
\item  
One of the relatively minimal models over $S/\langle 
\overline{g}_{3}, \overline{h} \rangle$ of 
$$\kappa_{3,1} : X_{3,1} \overset\Phi_{3,1}\to\longrightarrow 
\overline{X}_{3,1} = A_{3}/\langle g_{3}, h \rangle 
\overset\overline{\kappa}\to\longrightarrow S/\langle 
\overline{g}_{3}, \overline{h} \rangle,$$ 
where $\Phi_{3,1} : X_{3,1} \rightarrow \overline{X}_{3,1}$ 
is the contraction in (3.4)(2-1) and $\kappa$ is the morphism 
induced by the quotient map $A_{3} \rightarrow S := A_{3}/E$ given by the 
identity component $E$ of $\text{Ker}(h_{0} - id : A_{3} \rightarrow A_{3})$ 
of the Lie part $h_{0}$ of $h$ and $\overline{g}_{3}$ and $\overline{h}$ 
are the automorphisms of $S$ induced by $g_{3}$ and $h$. 
Moreover, in this case $I(W) = 3$ and $\rho(W) = 2$. 
\item 
One of the relatively minimal models over $S/G$ of 
$$p_{1} : \widetilde{(S \times E)/G} \overset\nu\to\rightarrow 
(S\times E)/G \overset\overline{p}_{1}\to \rightarrow S/G,$$ 
where $S$ is a normal K3 surface, $E$ is an elliptic curve, 
$G$ is a finite Gorenstein automorphism group of $S \times E$ 
whose element is of the form $(g_{S}, g_{E}) \in \text{Aut}(S) 
\times \text{Aut}(E)$ and $\nu$ is a crepant resolution 
of $(S\times E)/G$. Slightly more precisely, $G$ is of the form 
$G = H \rtimes \langle g \rangle$, where $H$ is a commutative group 
consisting of elements like $h = (h_{S}, t_{*})$
such that $\text{ord}(h_{S}) = \text{ord}(t_{*}) = \text{ord}(h)$ and 
$g$ is an element of the form $(g_{S}, \zeta_{I}^{-1})$ such that  
$g_{S}^{*}\omega_{S} = \zeta_{I}\omega_{S}$, where $I$ is the global 
canonical index of the base space $W$. Moreover $I \in \{2, 3, 4, 6 \}$.   
\endroster
\endproclaim
\remark{Remark 1} The cases (1), (2) and (3) do not appear 
in the classification 
in [Og4]. Indeed, $\pi_{1}(X)$ is infinite in the cases (1) and (2) and 
$\pi_{1}(X) \simeq C_{3}$ in the case (3).
\endremark 
\remark{Remark 2} In the case (4), the minimal resolution $S' \rightarrow S$ 
induces a birational morphism $(S' \times E)/G \rightarrow (S \times E)/G$. 
Since $(S' \times E)/G$ has only Gorenstein quotient singularities, 
$(S' \times E)/G$, hence $(S \times E)/G$, admits a crepant resolution 
[Ro, Main Result]. 
Moreover, it is well known that any two three-dimensional birational 
minimal models are connected by a sequence of flop and that three-dimensional 
flop does not affect singularities of minimal models [Kaw3] and [Kol]. 
Therefore, $X$ is smooth even in the cases (3) and (4).
\endremark 
\demo{Proof} As in [Og4, Section 2], let us consider the following 
commutative diagram: 
$$\CD 
X @<<< X \times_{W}T @<\nu<< Y\\
@V{\Phi}VV @V{\varphi'}VV @VV{\varphi_{Y}}V\\
W @<\pi<< T @= T
\endCD$$
where $\nu$ is a resolution of singularities of $X \times_{W}T$. 
Recall from [Og4, Pages 434-435] that $\text{Sing}(X \times_{W}T)$ 
is supported only in fibers of $\pi \circ \varphi'$. Let 
$\{E_{j} \vert j \in J\}$ 
be the set the two dimensional irreducible components in the fibers of 
$\varphi_{Y}$. Then $K_{Y} + \epsilon \sum_{j \in J} E_{j}$ is klt for small 
$\epsilon >0$ [Og4, Claim(2.10)] and we may therefore run the log Minimal 
Model Program for $Y$ with respect to $K_{Y} + \epsilon \sum_{j \in J} E_{j}$. 
By repeating the same procedure as in [Og4, Proposition (2.2)], we get, 
as its output, a contraction $\varphi_{Z} : Z \rightarrow T$ such that  
\roster
\item    
$\Cal O_{Z}(K_{Z}) \simeq \Cal O_{Z}$ and $h^{1}(\Cal O_{Z}) = 1$ [Og4, Pages 
435 - 436 (1) and (6)];
\item 
$\text{Sing}(Z)$ is purely one dimensional compound Du Val singularities 
along some fibers and is also equi-singular along each of such fibers 
(and therefore, $Z$ adimits at worst Gorenstein quotient singularities) 
[Og4, Lemma (2.7) and its proof]; 
\item 
the natural birational action of the Galois group 
$\langle g \rangle \simeq C_{I}$ on $Z$ is a biregular Gorenstein action and 
induces a birational map $\alpha : Z/\langle g \rangle \cdots \rightarrow X$ 
over $W$ [Og4, Page 436 (8) and (9)]. 
\endroster
{\it From now on, our argument differs to the one in [Og4].} By (2), $Z$ 
admits a unique crepant resolution $\mu : V \rightarrow Z$ 
(cf. [Re, Main Theorem (II)]). 
This $V$ is a smooth threefold such that 
$\Cal O_{V}(K_{V}) \simeq \Cal O_{V}$ and $h^{1}(\Cal O_{V}) = 1$.  
Moreover, the biregular action $\langle g \rangle \rightarrow \text{Aut}(Z)$ 
lifts to the one on $V$ again biregularly, because of the uniqueness of $V$. 
Denote by $a_{V} : V \rightarrow A$ the Albanese map of $V$ and 
set $\varphi_{V} := \varphi_{Z}\circ\mu : 
V \rightarrow T$. Then $A$ is an elliptic curve 
and $a_{V}$ is an \'etale fiber bundle over $A$ by [Kaw2, Theorem 8.3]. 
Moreover, by the uniqueness of the Albanese map, the action 
$\langle g \rangle \rightarrow \text{Aut}(V)$ descends equivariantly 
to the one on $A$. These actions make both $a_{V} : V \rightarrow A$ and 
$\varphi_{V} : V \rightarrow T$ $\langle g \rangle$-stable. 
Note that $\pi_{1}(V)$ is infinite by $h^{1}(\Cal O_{V}) > 0$.  
Therefore, by (2.1), $V$ admits the minimal 
splitting covering $\gamma : U \rightarrow V$  
from either an abelian threefold or the product of a K3 surface and 
an elliptic curve. Denote by $H$ the Galois group of $\gamma$. 
Since $\Cal O_{V}(K_{V}) \simeq \Cal O_{V}$, this $H$ is a Gorenstein 
automorphism group of $U$. 
Let us take the Stein factorisations of $a_{V} \circ \gamma$ and $\varphi_{V} 
\circ \gamma$ and denote them by: 
$$\CD 
S @< \varphi_U << U @> a_U >> E\\
@V \gamma_S VV @V \gamma VV @VV \gamma_E V\\
T @< \varphi_V << V @> a_V >> A.
\endCD$$
Then, by the uniqueness of the minimal splitting covering and by the 
uniqueness of the Stein factorisation, 
the action $\langle g \rangle \rightarrow \text{Aut}(V)$ again lifts to the 
one $\langle g \rangle \rightarrow \text{Aut}(U)$ equivariantly 
and also descends to both $\langle g \rangle \rightarrow \text{Aut}(S)$ 
and $\langle g \rangle \rightarrow \text{Aut}(E)$ and makes the diagram 
above $\langle g \rangle$-equivariant. Note that the action of $g$ on $U$ 
is also Gorensteion, because $\gamma$ is \'etale, and that the order of $g$ 
as an element of $\text{Aut}(U)$ is still $I$, because $\langle g \rangle 
\rightarrow \text{Aut}(V)$, and hence $\langle g \rangle 
\rightarrow \text{Aut}(U)$, is faithful. Similarly, $H \rightarrow 
\text{Aut}(U)$ descends to $H \rightarrow \text{Aut}(S)$ and 
$H \rightarrow \text{Aut}(E)$ to make the above diagram $H$-equivariant 
(where we define the actions of $H$ on the varieties in the bottom line as 
the trivial action). Moreover, by the connectedness of fibers of 
$\varphi_{V}$ and 
$a$ and by the finiteness of $\gamma$, $\gamma_{S}$, $\gamma_{E}$, 
we see that 
$\varphi_{V} : V \rightarrow T$ and $a_{V} : V \rightarrow A$ 
are isomorphic to the induced morphisms 
$\overline{\varphi}_{U} : U/V \rightarrow S/H$ and 
$\overline{a}_{U} : U/H \rightarrow E/H$ respectively. 
Set $G := \langle H, g \rangle$ as a subgroup of $\text{Aut}(U)$ 
and denote by $\rho_{G,S} : G \rightarrow \text{Aut}(S)$ 
and $\rho_{G,E} : G \rightarrow \text{Aut}(E)$ the equivariant actions 
found above. (Note that in a priori there is no reason why 
$\rho_{G,S}$ and $\rho_{G,E}$ are faithful.) 
\par
\vskip 4pt
Let us first treat the case where $U = S' \times E'$, the product of 
a K3 surface $S'$ and an ellptic curve $E'$. Recall by (2.25) 
that each element 
$h$ of $G$ is of the form $(h_{S'}, h_{E'}) \in \text{Aut}(S') \times 
\text{Aut}(E')$. In particular, the natural projections 
$p_{1} : U \rightarrow S'$ and $p_{2} : U \rightarrow 
E'$ are $G$-stable. We show that this case falls into the case (4) 
of (3.7).    
\proclaim{Claim (3.8)} 
\roster
\item 
The contractions $a_{U} : U \rightarrow E$ and $p_{2} : U \rightarrow E'$ 
are identically isomorphic. In particular, $E$ is an elliptic curve.  
\item 
The contraction $\varphi_{U} : U \rightarrow S$ factors through 
$p_{1}$, or more precisely, there exists a birational 
morphism $\tau : S' \rightarrow S$ such that $\varphi_{U} = \tau \circ p_{1}$.
\endroster 
In particular, the $G$-stable morphism 
$\varphi_{U} \times a_{U} : U \rightarrow S \times E$ is a crepant 
birational morphism. 
\endproclaim
\demo{Proof} Since $h^{1}(\Cal O_{S}) = 0$, we have 
$\text{Pic}(S' \times E') = 
p_{1}^{*}\text{Pic}(S') \otimes p_{2}^{*}\text{Pic}(E')$. Therefore, 
any divisor on $X$  
is linearly equivalent to a divisor of the form 
$D := p_{1}^{*}C + p_{2}^{*}L$, where $C$ and $L$ are divisors on $S'$ and 
$E'$ respectively. Note also that
$D$ is nef if and only if both $C$ and $L$ are nef. 
The morphism $a_{E}$ is given by such a nef divisor $D$ that $\nu(X, D) = 1$, 
because $\text{dim}(E) = 1$. Here and in what follows, we denote by 
$\nu(X, D)$ the numerical Kodaira dimension.
Therefore, $(\nu(S', C), \nu(E', L))$ is either $(0, 1)$ or $(1, 0)$, 
where $C$ and $L$ are same as above. 
However, in the latter case we have $\Phi_{D} = \Phi_{C} \circ p_{2}$ 
and the base space must then be $\Bbb P^{1}$, a contradiction. 
Hence, $(\nu(S', C), \nu(E', L)) = (0, 1)$ and $a_{U} = \Phi_{D}$ factors 
through $p_{2}$. Since both $a_{U}$ and $p_{2}$ have connected fibers, we 
get the assertion (1). 
Similarly, by $\text{dim}(S) = 2$, the contraction $\varphi_{S}$ is 
given by a nef divisor $D$ whose numerical Kodaira dimension is two. 
Therefore, $(\nu(S', C), \nu(E', L))$ is either $(2, 0)$ or $(1, 1)$ 
in this case. 
However, in the latter case we would have $S \simeq \Bbb P^{1} \times E'$, 
a contradiction. Hence, the first case occurs and therefore
$\varphi_{U} = \Phi_{D}$ factors through 
$p_{1}$. Again, since both $\varphi_{U}$ and $p_{1}$ have connected fibers, 
the induced morphism $S' \rightarrow S$ is birational. Now the last assertion 
follows from (1) and (2). \qed 
\enddemo
By (3.8) and (2.27), the restrictions $\rho_{G, E} \vert H$ and 
$\rho_{G, S} \vert H$ are both 
injective. Moreover, $\rho_{G, E}(H)$ is a translation subgroup of $E$, 
because $h^{0}(\Omega_{V}^{1}) = h^{0}(\Omega_{U}^{1}) = 1$. 
In particular, 
$\rho_{G, E}(H)$ is isomorphic to $C_{n} \oplus C_{m}$ for some 
$1 \leq n \vert m$.   
 Since $\langle g \rangle \rightarrow \text{Aut}(S)$ is also 
a lift of the original 
$C_{I} \simeq \langle g \rangle \hookrightarrow \text{Aut}(T)$ 
by the equivariantness, $\rho_{G, S} \vert \langle g \rangle$ is also 
injective. On the other hand, since both $S$ and $T$ are normal K3 surfaces, 
we have $\omega_{S} = \gamma_{S}^{*}\omega_{T}$. Therefore, by equivariantness 
and by $g^{*}\omega_{T} = \zeta_{I}\omega_{T}$, we have  
$g^{*}\omega_{S} = \zeta_{I}\omega_{S}$. Recall that $g^{*}\omega_{U} 
= \omega_{U}$. Then, by (3.8), 
$g^{*}\omega_{S \times E} = \omega_{S \times E}$. 
Therefore, $g^{*}\omega_{E} = \zeta_{I}^{-1}\omega_{E}$. This implies that
$\langle g \rangle \rightarrow \text{Aut}(E)$ is also injective and that
the image of $g$ is a Lie automorphism of $E$ of order $I$ under 
appropriate origin of $E$. 
Combining these together with the structure of automorphism group 
of an elliptic curve, we obtain 
$I \in \{2, 3, 4, 6\}$, 
a semi-direct decomposition $\rho_{G, E}(G) = \rho_{G, E}(H) \rtimes 
\rho_{G, E}(\langle g \rangle)$ and the injectivity 
$\rho_{G, E} : G \rightarrow 
\text{Aut}(E)$. This also implies $G = H \rtimes \langle g \rangle$ and 
gives an isomorphism $H \simeq C_{n} \oplus C_{m}$. 
Let us take $\tau \in \text{Ker}(\rho_{G, S})$ 
and write $\tau = h \circ g^{i}$ ($h \in H$). 
Since $\tau^{*}\omega_{S} = \omega_{S}$, we have $I \vert i$ and 
hence, $\tau = h$. Therefore $\tau = 1$ by the injectivity of 
$\rho_{G, S} \vert H$. This shows the injectivity of $\rho_{G, S}$. 
Now combining $G = H \rtimes \langle g \rangle$ with the construction, we 
readily see that the induced morphism 
$\overline{\gamma_{S} \circ \pi} : S/G \rightarrow W$ is an isomorphism 
and that the original $\Phi : X \rightarrow W$ and the induced contraction 
$\overline{p}_{1} : (S \times E)/G \rightarrow 
S/G$ are birationally isomorphic through the isomorphism 
$\overline{\gamma_{S} \circ \pi}$ and the composition of birational maps,
$\overline{\mu \circ \gamma \circ (\varphi_{U} \times a_{U})^{-1}} 
: (S \times E)/G \cdots\rightarrow Z/\langle g \rangle$ and $\alpha : 
Z/\langle g \rangle \cdots\rightarrow X$. 
Hence, $\Phi : X \rightarrow W$ falls in the Case (4) in (3.7). 
\par
\vskip 4pt   
Next we treat the case where $U$ is an abelian threefold. 
Our goal is to show that in this case $\Phi : X \rightarrow W$ falls into 
either one of cases (1), (2), (3) of (3.7). 
\par
By (3.1), $S$ is an abelian surface and $E$ is an elliptic curve. 
Note that $\varphi_{V}$ does not factor through $a_{V}$, because 
$h^{1}(\Cal O_{T}) = 0$. 
Then, $G$-stable map $\varphi_{U} \times a_{U} : U \rightarrow S 
\times E$ is surjective, and therefore, is an isogeny.
This, in particular, implies $g^{*}\omega_{E} = \zeta_{I}^{-1}\omega_{E}$, 
because $g^{*}\omega_{V} = \omega_{V}$ and $g^{*}\omega_{S} = 
\zeta_{I}\omega_{S}$. Therefore, $\rho_{G, E} \vert \langle g \rangle$ 
is injective.    
Assume that $\rho_{G, E} \vert H : H \rightarrow \text{Aut}(E)$ 
is not injective.  Then, there exists $1 \not= h \in \text{Ker}(\rho_{G, E} 
\vert H)$. Let $F$ be a fiber of $a_{U}$. Then $F$ is an abelian surface and 
$h$ acts on $F$. Since $h^{*}\omega_{U} = \omega_{U}$, we have 
$h^{*}\omega_{F} = \omega_{F}$. Moreover, $h$ has no fixed points.  
Therefore, $h \vert F$ must be a translation by the classification 
of automorphisms of abelian surfaces (eg. [Kat]). Thus,  
$h^{*} \vert H^{0}(U, \Omega_{U}^{1}) = id$ and $h$ must be also a 
translation of $U$. However, this contradicts 
the fact that $U \rightarrow V$ is the minimal splitting 
covering. Therefore, $\rho_{G,E} \vert H$ is injective. 
Note also that $\rho_{G, E}(H)$ is a translation group of $E$, 
because both $A = E/H$ and $A$ are elliptic curves. 
Now we may repeat the same argument as in the previous case to 
conclude that $I \in \{2, 3, 4, 6\}$, $\rho_{G, E}$ is injective, 
$H \simeq C_{n} \oplus C_{m}$ for some $1 \leq n \vert m$ 
and $G = H \rtimes \langle g \rangle$. Repeating the same argument for 
$\varphi_{U} : U \rightarrow S$, we also get the injectivity of $\rho_{G, S} 
\vert H$. Now the injectivity of $\rho_{G, S}$ again follows 
from the same argument as in the previous case. 
\proclaim{Claim (3.9)} $H \simeq C_{J}$, where $J \in  \{2, 3, 4, 6 \}$.
\endproclaim 
\demo{Proof} Put $H = \langle h_{1} \rangle \oplus \langle h_{2} \rangle 
(\simeq C_{n} \oplus C_{m})$. Since $H$ is Gorenstein and acts on $E$ 
as a translation, there exists a basis of $H^{0}(U, \Omega_{U}^{1})$ such that 
the matrix representation of $H$ is of the form, 
$h_{1}^{*} = \text{diag}(1, \zeta_{n}, \zeta_{n}^{-1})$ 
and 
$h_{2}^{*} = \text{diag}(1, \zeta_{m}, \zeta_{m}^{-1})$ 
(by changing the generaters if necessary). 
This implies $n, m \in \{1, 2, 3, 4, 6 \}$ as in (2.4)(3), 
and also $(h_{1} \circ h_{2}^{-m/n})^{*} = id$. Therefore, 
$h_{1} = h_{2}^{-m/n}$ by the definition 
of the minimal splitting covering. 
Thus $n = 1$, and hence, $H \simeq C_{m}$  
Since $h^{1}(\Cal O_{U}) = 3$ and $h^{1}(\Cal O_{V}) = 1$, 
we also see that $m \not= 1$. \qed
\enddemo
\enddemo
From now we argue dividing into cases according to the value $J$ in (3.9). 
Set $H = \langle h \rangle$. The basic idea of proof is to play ``fixed point 
game''. For proof, we also recall here that 
$\text{ord}(g) = I \in \{ 2, 3, 4, 6 \}$. 
\proclaim{Claim (3.10)} $J \not= 6$. 
\endproclaim 
\demo{Proof} Assume to the contrary, that $J = 6$. Take the origin $0$ of 
$E$ in $E^{g}$ and choose a global coordinate $z$ around $0$ of $E$. Then, 
we have $g(z) = \zeta_{I}^{-1}z$, $h(z) = z + p$, 
and $g^{-1}hg(z) = z + \zeta_{I}p$, where $p$ is a torsion point 
of oder $6$. In addition, since $\langle h \rangle$ is a normal 
group of $G$, we have 
either $g^{-1}hg = h$ or $g^{-1}hg = h^{-1}$. 
According to these two cases, we have $\zeta_{I}p = p$ and 
$-\zeta_{I}p = p$ respectively. Note that $E^{\zeta_{6}} =  
E^{-\zeta_{3}} = \{0\}$, 
$E^{-\zeta_{6}} = E^{\zeta_{3}}  \subset (E)_{3}$, $E^{\zeta_{4}} = 
E^{-\zeta_{4}} \subset E^{-1} = (E)_{2}$. 
Then, $I = 2$ and $g^{-1}hg = h^{-1}$, because $p$ is a point of order $6$.  
Let us consider the action of $G$ on $S$. Then $h^{*} 
\vert H^{0}(S, \Omega_{S}^{1})$ is of the form $\text{diag}(\zeta_{6}, 
\zeta_{6}^{-1})$ by the injectivity of the action 
an by $h^{*}\omega_{S} = \omega_{S}$. By this description 
and by the topological Lefschetz fixed point formula, we also see 
that $S^{h}$ is a one point set. Set $S^{h} = \{Q\}$. Then, $g(Q) = Q$ 
and there exist global coordinates $(x_{Q}, y_{Q})$ around $Q$ 
such that the (co-)action of $g$ is written as 
$g(x_{Q}, y_{Q}) = (-x_{Q}, y_{Q})$, because $g$ is an involution with 
$g^{*}\omega_{S} = -\omega_{S}$. Therefore, $(x_{Q} = 0)$ 
is a fixed curve of $g$. However, this contradicts 
the fact that $W = S/G$. Indeed, the quotient map $S \rightarrow W$ has no 
ramification curves, because $K_{S} \equiv 0$ and $K_{W} \equiv 0$. 
\qed   
\enddemo
\proclaim{Claim (3.11)} Assume that $J = 4$. Then $\Phi : X \rightarrow W$ 
falls into the case (2). 
\endproclaim 
\demo{Proof} As in (3.10), we may write the actions of $g$ and $h$ on $E$ 
as $g(z) = \zeta_{I}^{-1}z$ and $h(z) = z + p$, where $p \in E$ is a torsion 
point of order $4$. Then, by the same argument as the first part 
of the proof of (3.10), we see that $I = 2$ and $g^{-1}hg = h^{-1}$. 
In particular, $\langle h \rangle \rtimes \langle g \rangle \simeq D_{8}$. 
Let us consider the action of $G$ on $S$. Note that 
$\text{ord}(g^{i}h) = 2$ and 
$(gh^{i})^{*}\omega_{S} = -\omega_{S}$. Then, by the same argument as in 
the last part of the proof (3.10), we see that $g^{i}h$ has no fixed points 
on $S$. In particular, $g^{i}h$ has no fixed points on $U$. Note also 
that $h$ has also no fixed points on $U$, because $h$ is fixed point free on 
$E$. Therefore, the action of $G$ on $U$ has no fixed points and $U/G$ is then 
a Calabi-Yau threefold of Type A (0.1) (I) (2). 
Moreover, since $U/G$ is birational to $X$ 
and contains no rational curves, $U/G$ is isomorphic to $X$.
Now the rest of assertion follows from  Remark (2.22). 
\qed   
\enddemo
\proclaim{Claim (3.12)} Assume that $J = 3$. Then $\Phi : X \rightarrow W$ 
falls into the case (3). 
\endproclaim 
\demo{Proof} Again same as before, we may write the actions of 
$g$ and $h$ on $E$ 
as $g(z) = \zeta_{I}^{-1}z$ and $h(z) = z + p$, where $p \in (E)_{3} - \{0\}$. 
Then, again by the same argument as the first part 
of the proof of (3.10), we see that either $I = 2, 6$ and 
$g^{-1}hg = h^{-1}$, or $I = 3$ and $g^{-1}hg = h$ hold. 
Note that $h^{*} \vert H^{0}(S, \Omega_{S}^{1}) = 
\text{diag}(\zeta_{3}, \zeta_{3}^{2})$ under an appropriate basis 
$\langle v_{1}, v_{2} \rangle$. In particular, $\vert S^{h} \vert = 9$ by 
the Lefschetz fixed point formula.
First consider the case  $I = 2$ or $6$ and $g^{-1}hg = h^{-1}$. 
In this case, $g^{3}$ has a fixed point in $S^{h}$, 
because $\text{ord}(g^{3}) = 2$. Then, as in (3.10), 
$g^{3}$ have a fixed curve and gives the same contradiction. 
Therefore, $I = 3$ and $g^{-1}hg = h$, that is, $G = \langle h \rangle \oplus 
\langle g \rangle \simeq C_{3}^{\oplus 2}$ and $g^{*} 
\vert H^{0}(S, \Omega_{S}^{1}) = \text{diag}(\zeta_{3}^{2}, \zeta_{3}^{2})$ 
or $\text{diag}(1, \zeta_{3})$ under the same basis 
$\langle v_{1}, v_{2} \rangle$. By replacing $g$ by 
$gh^{2}$ if necessary, we may assume from the first that 
$g^{*} 
\vert H^{0}(S, \Omega_{S}^{1}) = \text{diag}(\zeta_{3}^{2}, \zeta_{3}^{2})$.  
Recall that $g^{*}\omega_{U} = \omega_{U}$ 
and $h^{*}\omega_{U} = \omega_{U}$. Then, $g$ acts on 
$U$ as a scalar multiplication by $\zeta_{3}^{2}$ and  
$h^{*} \vert H^{0}(U, \Omega_{U}^{1})$ is of the form 
$\text{diag}(1, \zeta_{3}, \zeta_{3}^{2})$. In particular, 
$U \simeq A_{3}$ by (3.2)(1). 
It remains to observe that the induced action of $h$ on 
$A_{3}/\langle g \rangle$ has no fixed points. Assume to the contrary, 
that $h$ has a fixed point $\overline{Q}$ on $A_{3}/\langle g \rangle$. 
Then there exists a point $Q \in A_{3}$ such that 
$g^{i}h(Q) = Q$ for some $i \in \{0, 1, 2\}$. 
Then $\varphi_{U}(Q) \in S$ is also a fixed point of $g^{i}h$. 
Here, we have $i \not= 0$, because $h$ has no fixed points on $E$ and hence on 
$U$. However, then the action $g^{i}h$ on $S$ has a fixed curve passing 
through $\varphi_{U}(Q)$, because $g^{i}h \vert H^{0}(S, \Omega_{S}^{1})$ 
has an eigen value $1$ if  $i \in \{ 1, 2 \}$, a contradiction. 
Therefore $h$ has no fixed points on  $A_{3}/\langle g \rangle$. 
\qed    
\enddemo
\proclaim{Claim (3.13)} Assume that $J = 2$. Then $\Phi : X \rightarrow W$ 
falls into the case (1). 
\endproclaim 
\demo{Proof} In this case $h^{-1} = h$. So, in a priori, $g^{-1}hg = h$. 
Again, aplying the same argument as the first part of the proof (3.10), 
we see that either $I = 2$ or $4$. Let us consider first the case 
where $I = 4$. Then, there exists 
a basis of $H^{0}(S, \Omega_{S}^{1})$ under which 
$h^{*} = \text{diag}(-1, -1)$ and $g^{*} = \text{diag}(1, \zeta_{4})$ 
or $\text{diag}(-1, -\zeta_{4})$. Replacing $g$ by $gh$ in the first case, 
we may assume from the first that $g^{*} = \text{diag}(-1, -\zeta_{4})$. 
Then $g$ has a fixed point on $S$ and therefore so does $g^{2}$. However, 
$g^{2}$ has then a fixed curve, a contradiction. Therefore $I = 2$ 
and $G \simeq C_{2}^{\oplus 2}$. In this case neither $g$ nor $gh$ has 
a fixed points on $S$, because they are non-Gorenstein involution on $S$. 
In addition, $h$ has no fixed points on $U$. 
Therefore $G$ acts freely on $U$. Now the same argument 
as the last part of (3.11) gives the result. 
\qed 
\enddemo
Now we are done. Q.E.D. of (3.7).
\head
{\S 4. Finiteness of $c_{2}$-contractions of a Calabi-Yau threefold} 
\endhead
\par 
\vskip 4pt 
In this section, we prove Theorem (0.4) in Introduction. For proof, 
it is convenient to introduce the notion of the maximal $c_{2}$-contraction: 

\proclaim{Lemma-Definition (4.1)}
There exists a $c_{2}$-contraction 
$\varphi_{0} : X \rightarrow W_{0}$ such that every $c_{2}$-contraction 
$\Phi : X \rightarrow W$ of $X$ factors through $\varphi_{0}$, that is, 
there is a morphism $\mu : W_{0} \rightarrow W$ such that 
$\Phi = \mu \circ \varphi_{0}$. 
Moreover, such $\varphi_{0} : X \rightarrow W_{0}$ is unique 
up to identical isomorphism. We call this $\varphi_{0} : X \rightarrow W_{0}$ 
the maximal $c_{2}$-contraction of $X$.
\endproclaim 
\demo{Proof} Let us choose a $c_{2}$-contraction 
$\varphi_{0} : X \rightarrow W_{0}$ 
such that $\rho(W_{0})$ is maximal among all $c_{2}$-contractions of $X$. 
We show that this $\varphi_{0}$ is the desired one. 
Let $\Phi : X \rightarrow W$ be any $c_{2}$-contraction. 
Take divisors $D_{0}$ and $D$ such that $\varphi_{0} = \Phi_{D_{0}}$ and 
$\Phi = \Phi_{D}$. Let us 
consider the contraction given by 
$m(D_{0} + D)$ for suitably large $m > 0$ and denote this contraction 
by $\Phi' := \Phi_{m(D_{0}+D)} : X \rightarrow W'$.  
Since $(c_{2}(X).D_{0}+D) =(c_{2}(X).D_{0}) + (c_{2}(X).D) = 0$, 
we see that $\Phi'$ is also a $c_{2}$-contraction. 
Moreover, by the construction, this $\Phi'$ factors through both 
$\Phi$ and $\varphi_{0}$, that is, there exist morphisms 
$p_{0} : W' \rightarrow W_{0}$ and 
$p : W' \rightarrow W$ such that 
$\Phi = p \circ \Phi'$ and $\varphi_{0} = p_{0} \circ \Phi'$. 
Note also that $W_{0}$ and $W'$ are both $\Bbb Q$-factorial 
by the classification of $c_{2}$-contractions. Indeed, 
they have at most quotient singularities (See Section 3).  
Hence, by the maximality of $\rho(W_{0})$, 
the morphism $p_{0}$ must be an isomorphism.  
Then, $\mu := p \circ p_{0}^{-1} : W_{0} \rightarrow W$ 
gives a desired factorisation. This argument also implies 
the last assertion. \qed
\enddemo

We proceed our proof of (0.4) dividing into cases according to the 
structure of the maximal $c_{2}$-contraction.   
In apriori, there are six possible cases: 
\par
Case I. $\varphi_{0}$ is an isomorphism (3.3) (= (0.1));
\par 
Case B. $\varphi_{0}$ is a birational contraction but not isomorphism 
(3.4);
\par 
Case K. $\varphi_{0}$ is of Type IIK (3.7);
\par 
Case A. $\varphi_{0}$ is of Type IIA (3.6);
\par 
Case P. $\text{dim}(W_{0}) = 1$; and 
\par 
Case T. $\text{dim}(W_{0}) = 0$. 
\par \noindent
In Case I, the result follows from (0.1)(IV). 
In Case P, $\varphi_{0}$ is an abelian fibration over $\Bbb P^{1}$ 
and this is the only (non-trivial) $c_{2}$-contraction of $X$. 
Case T is nothing but the case where $X$ admits no (non-trivial) 
$c_{2}$-contractions. It remains to consider Cases B, K, A.  
    
\demo{Proof of (0.4) in Case B} In this case $\varphi_{0} : 
X \rightarrow W_{0}$ is isomorphic to one of the contractions given in (3.4).  
First we treat the case where $\varphi_{0} : X \rightarrow W_{0}$ is 
isomorphic to $\Phi_{7} : X_{7} \rightarrow \overline{X}_{7}$ (3.4)(1).
Recall that $\Phi_{7}$ is the unique crepant resolution of 
$\overline{X}_{7}$. Therefore, it is sufficient to show the following:
\proclaim{Lemma (4.2)} $\overline{X}_{7}$ admits no non-trivial 
contractions. 
\endproclaim 
\demo{Proof} Let $f : \overline{X}_{7} \rightarrow W$ be a contraction 
and consider the Stein factorisation of the map 
$f \circ q$, where $q : A_{7} \rightarrow \overline{X}_{7}$ the quotient map: 
$$\CD 
A_{7}  @> f' >> V \\
@V q VV @VV q' V \\
\overline{X}_{7} @>> f > W.
\endCD$$  
Then $V$ is an abelian variety (3.1). Moreover, 
by the uniqueness of the Stein factorisation,  
$\langle g_{7} \rangle$ acts on $V$ equivariantly 
with respect to $f'$ and the induced map 
$\overline{f}' : \overline{X}_{7}
\rightarrow V/\langle g_{7} \rangle$ 
coincides with $f : \overline{X}_{7} \rightarrow W$. 
In addition, the action $g_{7}$ on $V$ is of order $7$, 
has only isolated fixed points, and satisfies 
$(g_{7})^{*}\omega_{V} \not= \omega_{V}$ if $\text{dim}V < 3$, 
because $\overline{f}'\circ\Phi_{7}$ is also a $c_{2}$-contraction. 
However, since $\varphi(7) = 6$, there are no elliptic curves 
and no abelian surfaces which admit such an automorphism. 
Therefore $V$ is an abelian threefold and $f'$ must be an isomorphism. 
\qed
\enddemo 
Next consider the case where 
$\varphi_{0} : X \rightarrow W_{0}$ is isomorphic to
$\Phi_{3} : X_{3} \rightarrow \overline{X}_{3}$ in (3.4)(2-0). 
Since $\Phi_{3}$ is the unique crepant resolution,  
the same argument as the first half part of (4.2) 
reduces our proof to the finiteness of $g_{3}$-stable 
contractions of $A_{3}$ up to $g_{3}$-equivariant 
isomorphisms. Therefore, the result follows from (3.2)(2). 
\par
\vskip 4pt
Let us consider the case where $\varphi_{0} : X \rightarrow W_{0}$ is 
isomorphic to the unique crepant resolution 
$\Phi_{3,1} : X_{3,1} \rightarrow \overline{X}_{3,1}$ (3.4)(2-1). 
Again as before, it is sufficient to show the finiteness of contractions 
of $\overline{X}_{3,1}$. In this case, we can say more: 
\proclaim{Lemma (4.3)} The nef cone $\overline{\Cal A}(\overline{X}_{3,1})$ 
is a rational simplicial cone and $\overline{X}_{3,1}$ admits exactly 6 
different non-trivial contractions. 
\endproclaim 
\demo{Proof} Our proof is quite similar to the one for (0.1) (IV) 
and we give just a sketch. 
Let us consider the elliptic curves $E_{i}$ ($1 \leq i \leq 3$) given 
as the identity components of the kernel of the endmorphisms,
$\text{Ker}(h_{0}\circ g_{3}^{3-i} - id : A_{3} \rightarrow A_{3})$.
Let $q_{i} : A_{3} \rightarrow S_{i} := A_{3}/E_{i}$ be the quotient map. 
Then the action of $\langle h, g_{3} \rangle$ descends equivariantly to 
the one on $S_{i}$, which we denote by $\langle \overline{h}, 
\overline{g}_{3} \rangle$. Then again taking the quotient of $S_{i}$ 
by the identity component $K_{i}$ of the kernel 
of $\overline{h} \circ (\overline{g}_{3})^{3-i} - id$, we finally 
obtain three different abelian fibrations $A_{3} \rightarrow 
B_{i} := S_{i}/K_{i}$. 
Moreover, these fibrations are $\langle h, g_{3} \rangle$-stable 
and therefore induce three differnet abelian fibrations 
$\varphi_{i} : \overline{X}_{3,1} \rightarrow \Bbb P^{1}$. 
Since $\rho(\overline{X}_{3,1}) = 3$, the rest of the proof 
is same as in (0.1)(IV). 
\qed      
\enddemo 
Finally, we consider the case where $\varphi_{0} : X \rightarrow W_{0}$ 
is isomorphic to the crepant resolution 
$\Phi_{3,2} : X_{3,2} \rightarrow \overline{X}_{3,2}$ (3.4)(2-2). 
However, in this case $\overline{X}_{3,2}$ admits no non-trivial 
contractions, because $\rho(\overline{X}_{3,2}) = 1$, and we are done.  
\qed 
\enddemo

\demo{Proof of (0.4) in Case K} By the case 
assumption, $\varphi_{0} : X \rightarrow W_{0}$ is isomorphic 
to one of the contractions in (3.7). In the first three cases of (3.7), 
we have $\rho(W_{0}) \leq 2$ so that $W_{0}$ admits 
at most two contractions and we are done. Let us consider the last case in 
(3.7). This is the essential case. In this case, 
$\varphi_{0} : X \rightarrow W_{0}$ 
is isomorphic to one of the relatively minimal models of   
$p_{1} : \widetilde{(S \times E)/G} \rightarrow S/G$. 
We denote this model by 
$f_{0} : Y \rightarrow B_{0} := S/G$ and fix a birational map 
$\rho_{0} : Y \cdot \cdot \cdot \rightarrow (S \times E)/G$ over $B_{0}$. 
Let $y_{i} : Y_{i} \rightarrow B_{0}$ ($i = 1, 2, ..., I$) be the complete 
representatives of the set of the relatively minimal models of 
$f_{0} : Y \rightarrow B_{0}$ modulo 
isomorphism over $B_{0}$. These are finite in number by virtue of the result 
of Kawamata [Kaw5, Theorem 3.6]. Indeed, 
since $Y$ itself is a Calabi-Yau threefold, 
$f_{0} : Y \rightarrow B_{0}$ is, in particular, a Calabi-Yau fiber 
space. Let $s_{j} : S \rightarrow S_{j}$ ($j = 1, 2, ..., J$) be the complete 
representatives of the set of $G$-stable contractions of $S$ modulo 
$G$-equivariant isomorphism. These are also finite in number as was 
shown in (1.10) (= (0.5)). 
We denote by $b_{j} : B_{0} \rightarrow B_{j} := S_{j}/G$ 
the contraction induced by $s_{j}$. In order 
to complete the proof, it is enough to show the following:

\proclaim{Lemma (4.4)} 
Every $c_{2}$-contraction of $Y$ is isomorphic to 
either one of $b_{j} \circ y_{i} : Y_{i} \rightarrow B_{j}$.
\endproclaim
\demo{Proof} Let $f : Y \rightarrow B$ be a $c_{2}$-contraction 
and $b : B_{0} \rightarrow B$ 
the factorisation of $f$ and write $f = b \circ f_{0}$. 
Let us denote by $q : S \rightarrow B_{0} = S/G$ the quotient map 
and consider the Stein factorisation of the map $b \circ q$: 
$$\CD 
S  @> b' >> C \\
@V q VV @VV q' V \\
B_{0} @>> b > B.
\endCD$$
As in (4.2), by the uniqueness of the Stein factorisation, 
$G$ acts on $C$ equivariantly and makes 
$b' : S \rightarrow C$ a $G$-stable contraction. Moreover,
the induced morphism $\overline{b}' : S/G \rightarrow C/G$ 
coincides with $b : B_{0} \rightarrow B$. 
Choose $j \in \{1, 2, ..., J\}$ such that $b' : S \rightarrow C$ is 
$G$-equivariantly isomorphic to $s_{j} : S \rightarrow S_{j}$ and denote 
this isomorphism by: 
$$\CD 
S  @> \sigma >> S \\
@V b' VV @VV s_{j} V \\
C @>> \sigma_{C} > S_{j}.
\endCD$$ 
This pair $(\sigma, \sigma_{C})$ descends to a pair of isomorphisms, 
$\sigma_{B_{0}} : B_{0} \rightarrow B_{0}$ and  
$\sigma_{B} : B \rightarrow B_{j}$ which give an isomorphism 
between $b : B_{0} \rightarrow B$ and $b_{j} : B_{0} \rightarrow B_{j}$:
$$\CD 
B_{0}  @> \sigma_{B_{0}} >> B_{0} \\
@V b VV @VV b_{j} V \\
B @>> \sigma_{B} > B_{j}.
\endCD$$
Let us consider an automorphism  
$\tilde{\sigma} = (\sigma, id)$ of $S \times E$. 
By the description of the elements of $G$ (3.7)(4), the pair 
$(\tilde{\sigma}, \sigma)$ gives a $G$-equivariant isomorphism of 
$p_{1} : S \times E \rightarrow S$: 
$$\CD 
S \times E  @> \tilde{\sigma} >> S \times E \\
@V p_{1} VV @VV p_{1} V \\
S @>> \sigma > S.
\endCD$$
Therefore $\tilde{\sigma}$ induces an isomorphism, 
$\tau : (S \times E)/G \rightarrow (S \times E)/G$ 
such that the pair $(\tau, \sigma_{B_{0}})$ 
gives an isomorphism on $p_{1} : (S \times E)/G \rightarrow B_{0}$:
$$\CD 
(S \times E)/G  @> \tau >> (S \times E)/G \\
@V p_{1} VV @VV p_{1} V \\
B_{0} @>> \sigma_{B_{0}} > B_{0}.
\endCD$$
Note that $\sigma_{B_{0}} \circ f_{0} : Y \rightarrow B_{0}$ is a 
relatively minimal 
model of $f_{0} : Y \rightarrow B_{0}$ via 
the birational map $(\rho_{0})^{-1} \circ \tau \circ \rho_{0}$. 
Then, there exists 
$i \in \{1, 2, ... , I\}$ such that  
$\sigma_{B_{0}} \circ f_{0} : Y \rightarrow B_{0}$ and 
$y_{i} : Y_{i} \rightarrow B_{0}$ are isomorphic over $B_{0}$. 
Let us choose one of such isomorphisms and denote it by $\tau_{i} : Y 
\rightarrow Y_{i}$. Then, the pair 
$(\tau_{i}, \sigma_{B_{0}})$ gives an isomorphism 
between $f_{0} : Y \rightarrow B_{0}$ and $y_{i} : Y_{i} \rightarrow B_{0}$:
$$\CD 
Y  @> \tau_{i} >> Y_{i} \\
@V f_{0} VV @VV y_{i} V \\
B_{0} @>> \sigma_{B_{0}} > B_{0}.
\endCD$$
Composing $(\tau_{i}, \sigma_{B_{0}})$ with the pair 
$(\sigma_{B_{0}}, \sigma_{B})$, 
we get an isomorphism 
between $f = b \circ f_{0} : Y \rightarrow B$ and 
$b_{j} \circ y_{i} : Y_{i} \rightarrow B_{j}$. \qed
\enddemo 
This completes the proof in Case K. \qed 
\enddemo
\demo{Proof of (0.4) in Case A}
Finally we consider the case where $\varphi_{0} : X \rightarrow W_{0}$ is 
of the form 
$f_{0} : Y \rightarrow B_{0} := E_{\zeta_{3}}^{2}/
\langle \text{diag}(\zeta_{3}, \zeta_{3}) \rangle$ 
described in (3.6). 
Recall by (3.6) that the number of the
relatively minimal models of 
$f_{0} : Y \rightarrow B_{0}$ is just $2^{9}$. 
We denote them by 
$y_{i} : Y_{i} \rightarrow B_{0}$ ($i = 1, 2, ..., 2^{9}$). 
Let $p_{1} : B_{0} \rightarrow E_{\zeta_{3}}/\langle \zeta_{3} \rangle 
= \Bbb P^{1}$ be the natural projection to the first factor. 
In order to conclude the result, it is sufficient to show the following: 
\proclaim{Lemma (4.5)} Every $c_{2}$-contraction of $Y$ is isomorphic to 
one of 
$y_{i} : Y_{i} \rightarrow B_{0}$ and 
$p_{1} \circ y_{i} : Y_{i} \rightarrow \Bbb P^{1}$.
\endproclaim
\demo{Proof} By (3.2)(2), we see 
that each $\langle \text{diag}(\zeta_{3}, \zeta_{3}) \rangle$-stable 
(non-trivial) contraction of $E_{\zeta_{3}}^{2}$ 
is $\langle \text{diag}(\zeta_{3}, \zeta_{3}) \rangle$-equivariantly 
isomorphic to either $id : E_{\zeta_{3}}^{2} \rightarrow 
E_{\zeta_{3}}^{2}$ or 
$p_{1} : E_{\zeta_{3}}^{2} \rightarrow E_{\zeta_{3}}$. 
Now we may repeat the same argument as in (4.4) to obtain the result. 
\qed
\enddemo
This completes the proof of (0.4). \qed. 
\enddemo

\enddocument